\documentclass{amsart}

\usepackage{amsmath}
\usepackage{amsfonts}
\usepackage{amsthm}

\usepackage{graphicx} 
\usepackage{amssymb,mathrsfs}

\usepackage{color}

\theoremstyle{definition}

\theoremstyle{remark}

\numberwithin{equation}{section}

\newcommand{\Var}{\operatorname{Var}}
\newcommand{\E}{\operatorname{E}}
\newcommand{\N}{\mathbb{N}}
\newcommand{\R}{\mathbb R}

\begin{document}


\title[Themes in the work of Manfred Denker]
{From Ergodic Theory and Probability to Fractal Geometry and Dynamics: 
Themes in the Work of Manfred Denker}

\author[S. Boyd]{Suzanne Boyd}
\address{Department of Mathematical Sciences\\
University of Wisconsin Milwaukee\\
PO Box 413\\
Milwaukee, WI 53201, 
USA}
\email{sboyd@uwm.edu, ORCID: 0000-0002-9480-4848}

\author[H. Dehling]{Herold Dehling}
\address{Fakult\"at f\"ur Mathematik\\
Ruhr-Universität Bochum\\
Universitätsstraße 150\\
44780 Bochum\\
Germany
}
\email{herold.dehling@ruhr-uni-bochum.de}

\author[M. Schmoll]{Martin Schmoll}
\address{School of Mathematical and Statistical Sciences\\
Clemson University\\
Clemson, SC 29634\\ 
USA}
\email{schmoll@clemson.edu}
\thanks{Schmoll was partially supported by grant SFI-MPS-TSM-00013668 from Simons Foundation International.}

\author[M. Stadlbauer]{Manuel Stadlbauer}
\address{Departamento de Matem\'{a}tica\\
Universidade Federal do Rio de Janeiro\\
Ilha do Fund\~{a}o\\
21941-909 Rio de Janeiro (RJ)\\
Brazil}
\email{manuel@im.ufrj.br, ORCID: 0000-0003-2537-9128}
\thanks{Stadlbauer was partially supported by the Fundação Carlos Chagas Filho de Amparo à Pesquisa do Estado do Rio de Janeiro (FAPERJ) through grant E-26/204.324/2024.}

\author[C. Wolf]{Christian Wolf}
\address{Department of Mathematics and Statistics\\
Mississippi State University\\
Starkville, MS 39759, USA}
\email{cwolf@math.msstate.edu, ORCID: 0000-0002-7976-3574}
\thanks{Wolf was partially supported by grant SFI-MPS-TSM-00013897 from Simons Foundation International.}



\date{\today}
\begin{abstract}
This article surveys the mathematical contributions of Manfred Denker, with a focus on themes that connect ergodic theory, probability theory, dynamical systems, fractal geometry, and statistics. Denker’s highly influential work includes a systematic study of the statistical properties of dynamical systems, the development of limit theorems for dependent processes, and the use of thermodynamic formalism to relate geometric and measure-theoretic properties.

Particular emphasis is placed on the emergence of probabilistic behavior in deterministic systems, including central limit theorems, invariance principles or local limit theorems, under weak dependence assumptions or in infinite measure.
Further topics include equilibrium states and transfer operator methods, the role of conformal measures in fractal geometry, and the asymptotic theory of statistical procedures for dependent data, such as rank statistics and $U$--statistics.

In addition to these theoretical developments, the survey highlights contributions connecting rigorous analysis with computational and statistical methods. Taken together, these works illustrate a unifying perspective in which ergodic, probabilistic, geometric, and statistical methods interact in the study of dynamical systems.
\end{abstract}

\maketitle

\section{Introduction}
\label{sec:intro}

The interaction between probability, geometry, and dynamics has been one of the defining themes of modern mathematics over the past half century. Within this broad development, the work of Manfred Denker occupies a particularly distinguished place. Through a long and influential career, Denker has made fundamental contributions to ergodic theory, probability and statistics, dynamical systems, and fractal geometry, while simultaneously helping to shape a conceptual framework that unifies these disciplines. His work exemplifies the modern viewpoint that deterministic systems may exhibit statistical behavior of remarkable regularity and universality, and that geometric complexity can often be understood through probabilistic and dynamical principles.

The present article is a survey of Denker’s mathematical contributions, with particular emphasis on the themes that have shaped his research over several decades. Rather than attempting a complete account, we focus on selected areas in which his work has had lasting influence, highlighting both foundational results and the methods that connect different parts of mathematics. The article is organized into sections reflecting these themes, including ergodic theory, probability theory, fractal geometry, dynamical systems,  statistics, and computational methods.

One of the central threads in Denker’s research concerns the statistical properties of dynamical systems. Building upon the foundations of ergodic theory and thermodynamic formalism, he developed influential results on invariant measures, mixing properties, Gibbs states, and limit theorems for dependent processes generated by deterministic dynamics. In collaboration with numerous coauthors, Denker helped establish central limit theorems, laws of the iterated logarithm, and large deviation principles in settings far removed from classical independent random variables. These works demonstrated that chaotic dynamical systems possess intrinsic stochastic features and provided rigorous methods for quantifying statistical behavior in deterministic environments.

Closely related to these developments are Denker’s contributions to probability theory and mathematical statistics. A recurring theme is the analysis of stochastic processes with strong dependence structures, especially those arising from dynamical systems. By adapting empirical process methods and asymptotic techniques to dependent observations, he significantly broadened the scope of modern statistics. His work clarified the interplay between mixing conditions, decay of correlations, and statistical inference, influencing both probability theory and statistical mechanics.

Another important direction in Denker’s work concerns dimension theory and fractal geometry. The emergence of fractal sets as natural objects in nonlinear dynamics led to fundamental questions about Hausdorff dimension, multifractal analysis, and invariant geometric structures. Denker’s contributions revealed deep connections between thermodynamic formalism and geometric measure theory, showing how entropy, pressure, and invariant measures govern the fine structure of dynamically defined sets. These ideas have had lasting impact not only within pure mathematics, but also in the study of physical systems exhibiting chaotic or self-similar behavior.

The breadth of Denker’s work is reflected not only in the range of topics he has addressed, but also in the diversity of methods he has developed. A distinctive feature of his research is the continuous interplay between abstract theory and concrete applications, including computational and statistical approaches to dynamical systems.

Throughout his career, Manfred Denker  has worked in close connection with a wide academic network. As an undergraduate, and later as a graduate student and assistant professor, all in Erlangen, he was part of the large research group of his Ph.D.\ advisor Konrad Jacobs (1928--2015). In the 1960s Jacobs was the towering figure in Ergodic Theory and Information Theory in Germany.  Among Jacobs'  former students are Rudolf Ahlswede, Thomas Beth, Ernst Eberlein, Hans F\"ollmer, Hans-Otto Georgii, Mike Keane, Ulrich Krengel,  Dietrich Werner M\"uller, Hans-J\"urgen Schuh, Volker Strassen, and Manfred Denker. 
Already during his Erlangen period (1966--1974), Denker began to build his own network. 
As an undergraduate student, he spent the academic year 1969-70 at Warwick University, where he connected to Rufus Bowen, Bill Parry, and Peter Walters. It was through their work and through discussions with them that he received inspiration for his Diplomarbeit and also for his Ph.D.\ thesis. 
Toward the  end of his time in Erlangen, Denker spent the academic year 1973-74 at the University of Rennes, where Mike Keane and Jean-Pierre Conze were his closest colleagues. 

After starting his professorship in G\"ottingen in 1974, Manfred Denker's research network expanded quickly. Many researchers came as guests to G\"ottingen, perhaps only for a talk in the weekly research seminar, or a stay lasting from a week to an entire summer, or in some cases even for a year or more.  Among the long-term guests were, e.g.,  Jon Aaronson, Rabi Bhattacharya, Bob Burton, Alby Fisher, Misha Gordin, Adam Jakubowski,  Mike Keane, Gerhard Keller, Yuri Kifer, Gerhard Knieper, Zbigniew Nitecki, Feliks Przytycki, Hiroshi Sato, Mariusz Urba\'nski, Mishiko Yuri, Simon Waddington, and Wojbor Woyczy\'nski. 
Financial support came mostly through grants from, e.g.,  the German Science Foundation (DFG), the Alexander von Humboldt Foundation, the German Academic Exchange Service (DAAD), the Japanese Society for the Promotion of Science, the Volkswagen Stiftung, and the European Union. The financial resources expanded widely during the time of the Sonderforschungsbereich SFB 170 {\em Geometry and Analysis}  (1984--1996), where Manfred Denker was one of the PIs, allowing him to employ postdocs and Ph.D.\ students, and invite large numbers of guests. Many joint papers were written as a result of cooperations  with visiting researchers, often also involving local Ph.D.\ students and postdocs. 

At the same time, Manfred Denker frequently visited foreign universities. He spent longer periods at Case Western Reserve University in Cleveland, where he worked mostly with Wojbor Woyczy\'nski (1943--2021), at the University of Illinois at Urbana-Champaign, where he collaborated with Walter Philipp (1936--2006), at Oregon State University in Corvallis, where he worked with Bob Burton (1948--2018),  at Indiana University in Bloomington, where he collaborated with Madan Puri,  at Kyushu University, where he worked with Hiroshi Sato, and at the Federal University of Rio de Janeiro, where he worked with Samuel Senti.  Denker's students and postdocs profited enormously from being embedded into his large network. Many of them later became professors at major research universities, thus contributing to a further expansion of Denker's network. After retiring from his G\"ottingen professorship in 2008, Manfred Denker moved for 10 years to The Pennsylvania State University, where he built a new network of young researchers around him.

Finally, we would like to emphasize the personal influence of Manfred Denker. To all the authors of this survey, he has been an inspiring mentor, collaborator, and friend, and for some of us, this relationship has spanned several decades. His enthusiasm for mathematics, his openness to new ideas, and his dedication to collaboration have shaped a large international community of researchers whose work continues to build on his contributions.
We are deeply grateful for the opportunity to have learned from and collaborated with Manfred Denker, and we feel fortunate to count him not only as a colleague and mentor, but also as a friend.

\section{Ergodic Theory}
\label{sec:ergodic_theory}

This survey begins by presenting  some of Manfred Denker’s most prominent results in ergodic theory, starting with his work on topological models for ergodic automorphisms which was part of his Ph.D.\ thesis in the early 1970's.

\subsection{Strict ergodicity}
In the early 1970s, Jewett and Krieger discovered that every ergodic automorphism can be realized, up to measure-theoretic conjugacy, as a uniquely ergodic and minimal homeomorphism of a compact space. In other words, every such system admits a strictly ergodic topological model. Moreover, their work revealed that ergodic automorphisms can be encoded symbolically through the construction of finite generating partitions.

These ideas were taken further in Denker’s 1972 Ph.D.\ thesis, where the generator method was unified and sharpened, which gave rise to short and transparent proofs of Krieger's generator theorem and the existence of strictly ergodic models (\cite{MR352402,MR382588,MR382587}).

It is then natural to ask whether the same holds for ergodic and measure-preserving flows, whose positive answer nowadays is known as the Denker-Eberlein-Theorem,
published as ``{Ergodic flows are strictly ergodic}" in  \emph{Advances in Mathematics} in 1974 \cite{MR352403}.

It is worth noting that strict ergodicity for flows attracted considerable attention in the early 1970s, largely due to the seminal 1972 result of Furstenberg. In that work, he provided the first natural example by showing that the horocycle flow on compact hyperbolic surfaces is strictly ergodic. From a measure-theoretic perspective, however, this phenomenon is far from exceptional as the Denker–Eberlein theorem shows that every ergodic flow admits a strictly ergodic model.

\subsection{Finitary isomorphisms}
The above shows that strict or unique ergodicity is not invariant under measure-theoretic conjugacy. It is therefore natural to seek a notion of equivalence under which these properties are preserved. In his thesis, Denker already considered isomorphisms which, on the symbolic level, depend only on finitely many symbols, although the corresponding range may depend on the point and need not be uniformly bounded. This concept, today known as finitary isomorphisms, came into sharper focus in the late 1970s, with Keane and Smorodinsky (\cite{MR528969}).

From a structural point of view, joint work with Michael Keane showed that strict ergodicity is invariant under finitary isomorphisms, and contains the remarkable observation that integrability of the pointwise coding range implies a central limit theorem. In addition, they proved that every ergodic dynamical system with finite entropy admits a finite and finitary generator, a result now known as the finitary generator theorem. These important results were published as ``{Almost topological dynamical systems}" in the \emph{Israel Journal of Mathematics} in 1979 \cite{MR571401}.

At this point, it is also worth mentioning that the concept of finitary isomorphisms was already employed in the 1976 monograph ``{Ergodic theory on compact spaces}" by Grillenberger, Denker and Sigmund from 1976 \cite{MR457675} to provide an elegant proof of Krieger's generator theorem, and that this book has since become one of the classic text books on ergodic theory.

\subsection{Infinite ergodic theory}
After the foundational ratio ergodic theorem by Hopf and its far reaching refinements by Chacon and Ornstein, 
it was natural to seek pointwise ergodic theorems in which the ratio limits are replaced by suitably normalized ergodic sums. However, as observed by Jon Aaronson during the late 1970s and 1980s, infinite measure-preserving systems exhibit fundamentally different qualitative behavior. In particular, as classical normalization schemes drastically fail, he established refined notions of ergodicity like rational ergodicity and pointwise dual ergodicity.

Probably due to their friendship and mutual interests in abstract ergodic theory and probability theory, Denker and Aaronson began a fruitful collaboration that led to several influential papers on ergodic theory, probability theory, and their interactions; such as, their work on Markov fibred systems and their applications to parabolic rational maps from 1993 \cite{MR1107025}, and on local limit theorems for partial sums generated by Gibbs-Markov maps from 2001 \cite{MR1840194}.

In the context of abstract ergodic theory, they considered a broad class of non-invertible transformations with infinite invariant measures, which contains various classical examples from number theory, topological Markov chains and parabolic rational maps. This class is defined through the presence of a Darling-Kac set, which by definition allows one to control the returns to the set uniformly in average through the transfer operator, and the additional assumption that the associated return sequence $(a(n))$ is slowly varying.

Under these natural assumptions, Aaronson and Denker were able to establish a  class theorem for the almost sure behavior of the $n$-th partial sum, normalized by some sequence $(\varphi(n))$ in the paper titled ``{Upper bounds for ergodic sums of infinite measure preserving transformations}" published in  \emph{Transactions of the American Mathematical Society} in 1990 \cite{MR1024766}.
Depending on a summability condition on $(\varphi(n))$, the limit  superior is then almost surely either bounded from above or below by the integral of the observable. In particular, for the limiting case, that is $\varphi(n)  =  a (n /\log \log n) \log \log n$, one obtains a law of the iterated logarithm in infinite measure, which is related to the invariance principle and the LIL for stable laws.

A further important result on infinite ergodic theory, this time joint with  Aaronson and Albert Fisher, provides a solution for the absence of almost sure convergence of partial sums in infinite measure. In the above setting, they obtained almost sure convergence of a second order limit in $\log$-average, that is of
\[
\lim_{m \to \infty} \frac{1}{\log m} \sum_{n=1}^{m} \frac{1}{na(n)} \sum_{k=1}^n f\circ T^k.
\]
This second order theorem represents a generalization of the Birkhoff ergodic theorem to infinite measure preserving systems and is the first of its kind. It was published as
``{Second order ergodic theorems for ergodic transformations of infinite measure spaces}" in the  \emph{Proceedings of the American Mathematical Society} in 1992 \cite{MR1099339}.

\subsection{Equilibrium states} The construction of invariant probability measures with strong statistical properties is a central problem in ergodic theory. Furthermore, the construction of equilibrium states is particularly prominent, as the variational principle links topology, geometry, and measure theory, while the transfer operator provides direct access to  dynamical properties.

This approach, also known as thermodynamic formalism due to the analogies with statistical physics, was  successfully applied in the 1970s to hyperbolic dynamical systems with finite Markov partitions. However, in order to understand systems with some weak hyperbolicity, like maps from the quadratic family or with intermittent behavior, many researchers put their emphasis on non-invertible, piecewise monotone dynamical systems defined on the unit interval. Under the assumption that $h_{\mathrm{top}} > \sup \varphi - \inf \varphi$, where $\varphi$ is a relevant potential function, Hofbauer und Keller showed the existence of an equilibrium state in \cite{MR684242}.

In collaboration with Mariusz Urba\'nski, Denker began to develop and extend these ideas to rational maps in a series of influential papers. After providing a very general construction of conformal measures (\cite{MR1014246}), they proved in ``{Ergodic theory of equilibrium states for rational maps} "(\emph{Nonlinearity} 1991, \cite{MR1389624}) that the weaker condition $P(\varphi) > \sup \varphi$ suffices to guarantee the existence of an equilibrium state, this time for analytic endomorphisms of the Riemann sphere instead of a piecewise monotonic map on the interval, and that the equilibrium state is exact.  This result later was refined by Denker in joint work with Feliks Przytycki and Mariusz Urbański by pointing out that the transfer operator in fact has a spectral gap, which was published as
``{On the transfer operator for rational functions on the Riemann sphere}" in \emph{Ergodic Theory and Dynamical Systems} in 1996 \cite{MR1389624}.

\subsection{A multivariate ergodic theorem.} Another important collaborator of Denker was Mikhail Gordin (1944--2015), who is well-known for his proof of the central limit theorem for dynamical systems through martingale approximation. Gordin spent multiple summers in Göttingen. In ``{Limit theorems for von Mises statistics of a measure preserving transformation}", published in {\it Probability Theory and Related Fields} in 2014 \cite{MR3256808} they established a multivariate ergodic theorem by showing that, under certain conditions, the von Mises statistics
\[
x \mapsto \sum_{0\le i_1,\ldots,i_d\le n}
f(T^{i_1}x,\ldots,T^{i_d}x) , \quad n\ge 1\]
satisfies the almost sure and statistical ergodic theorem  under a one-dimensional action.
The significance of this result lies in its study of the measure class of multidimensional functions that embed into
$L^2(\mu^d)$ and admit a $\mathbb Z^d$-action
via coordinate-wise Koopman operators. It fills a notable gap in the literature, as previous work focused on individual functions of this type rather than on their measure-theoretic classification within $L^2$. The result was originally drafted around 2004 and completed at Penn State in 2013. Around that time,  Gordin was diagnosed with cancer and, after two years, sadly passed away. 

At this point, it should be noted that the research focus of Denker and Gordin was on probabilistic aspects of dynamical systems, which led to two influential papers in this area \cite{MR1736956,MR2048211}.

\subsection{Perspective}
Manfred Denker’s contributions to ergodic theory span a broad range of directions. His first influential results concerned strict ergodicity and finitary isomorphisms, which are classical topics in ergodic theory. Thereafter, he studied measure-preserving dynamical systems using probabilistic techniques and the theory of Markov operators, leading to new almost sure limit laws and ergodic theorems, in both finite and infinite measure settings, as well as the construction of invariant Gibbs measures.

A common characteristic of this second phase is the intricate use of operator-theoretic methods, sometimes combined with renewal theory. This approach allows for the effective analysis of geometric, dynamical, and statistical properties across a wide variety of settings, including non-uniformly expanding and parabolic dynamical systems.

\section {Probability Theory}
\label{sec:probability}

When starting his professorship in G\"{o}ttingen in 1974, Manfred Denker began to  teach probability and statistics, which shifted his focus from abstract ergodic theory toward analytic methods and tools present in these fields of Mathematics. Denker's teaching has profoundly influenced generations of students in G\"ottingen, for example through the lecture notes of his 1976 course on probability theory. These lecture notes contained many advanced topics such as the ergodic theorem, martingale theory, basics of stochastic processes, stable and infinitely divisible laws. These lecture notes, known as the \emph{Gr\"{u}ne Buch}, remained in use for several decades, and  became well known even beyond G\"{o}ttingen.

On the other hand, this also broadened Manfred’s perspective on ergodic theory and dynamical systems, as he began to consider probabilistic aspects of such systems. Among other things, this led to the question of the minimal conditions required for obtaining a central limit theorem or other limit laws.

\subsection{Limit theorems for weakly dependent processes}
While the ergodic theorem holds for an arbitrary stationary stochastic processes $(X_j)_{j\geq 1}$, more refined limit theorems, such as the central limit theorem and the law of the iterated logarithm, require further assumptions on the degree of dependence of the process, classically expressed by conditions on mixing coefficients. The latter arise from measures of the degree of dependence of two $\sigma$-fields $\mathcal{A}$, $\mathcal{B}$, the best known of which are
\begin{align*}
\alpha(\mathcal{A},\mathcal{B})&= \sup_{A\in \mathcal{A}, B\in \mathcal{B} } |P(A\cap B)-P(A)\, P(B)|, \\
\phi(\mathcal{A},\mathcal{B})&= \sup_{A\in \mathcal{A}, B\in \mathcal{B}} |P(B|A)- P(B)|,\\
\psi(\mathcal{A},\mathcal{B})&= \sup_{A\in \mathcal{A}, B\in \mathcal{B}} \left|\tfrac{P(A\cap B)}{P(A)P(B)} - 1\right|,
\end{align*}
together with
\[
\beta(\mathcal{A},\mathcal{B})= \sup \sum_{i=1}^k \sum_{j=1}^l |P(A_i\cap B_j)-P(A_i)\, P(B_j)|,
\]
where the supremum in the definition of $\beta(\mathcal{A},\mathcal{B})$ above extends over all partitions $ A_1,\ldots,A_k \in \mathcal{A}$ and $B_1,\ldots,B_l \in \mathcal{B} $ of the probability space $\Omega$. These measures give rise to the mixing coefficients of a stochastic process $(X_i)_{i\geq 1}$ measuring the degree of dependence of the $\sigma$-fields $\mathcal{F}_k^l=\sigma(X_k,\ldots,X_l)$ for sets of indices separated by some large $n$. More specifically, we define the $\alpha$-mixing coefficient
\[
  \alpha(n)=\sup_{k} \alpha(\mathcal{F}_1^k, \mathcal{F}_{k+n}^\infty),
\]
and similarly the $\beta$-, $\phi$- and $\psi$-mixing coefficients $\beta(n)$,  $\phi(n)$ and $\psi(n)$. The process $(X_i)_{i \geq 1}$ is called $\alpha$-mixing if $\lim_{n\rightarrow \infty} \alpha(n)=0$, and analogously one defines $\beta$-, $\phi$ and $\psi$-mixing processes. Alternative notations are strong mixing (instead of  $\alpha$-mixing), uniform mixing (instead of $\phi$-mixing) and absolutely regular (instead of $\beta$-mixing). Note that $\psi$-mixing implies $\phi$-mixing, and $\phi$-mixing implies $\beta$-mixing, which again implies $\alpha$-mixing.
Central limit theorems usually require moment conditions together with conditions on the rate of decay of one of the mixing coefficients. For example, Ibragimov (1975) \cite{MR362448} proved that for $\phi$-mixing processes $(X_i)_{i\geq 1}$ with  $\mathrm{Var}(X_1)<\infty$, satisfying

\begin{align*}
  \sum_{n=1}^\infty \phi^{1/2}(2^n)   <\infty, \quad
   \sigma_n^2:=\Var\big(\sum_{i=1}^n X_i\big)  \longrightarrow \infty,
\end{align*}
the central limit theorem holds, i.e., that
\[
  \frac{1}{\sigma_n} \sum_{i=1}^n (X_i-\E(X_1)) \stackrel{\mathcal{D}}{\longrightarrow} N(0,1).
\]

In their book, Ibragimov and Linnik (1971) \cite{MR322926} formulated the conjecture that the CLT should hold for $\phi$-mixing processes satisfying $\sigma_n^2\rightarrow \infty$, without any further assumptions on the rate of convergence of the $\phi$-mixing coefficients. This has become known as Ibragimov's  conjecture, which still stands open as we complete this article.
In two papers written around 1980, Manfred Denker tried to find minimal conditions for the CLT for stationary stochastic processes. In the short note ``{Uniform integrability and the central limit theorem for strongly mixing processes}" (\cite{MR899993}, 1986), originally presented in 1978 in lectures on stationary processes, Denker showed that for  $\alpha$-mixing processes uniform integrability of $ {\sigma_n^{-2}}(\sum_{i=1}^n (X_i-\E(X_i)) )^2$ is a necessary and sufficient condition for the CLT to hold.
Herold Dehling, Manfred Denker and Walter Philipp  obtained in ``{Central limit theorems for mixing sequences of random variables under minimal
conditions}", published in \emph{Annals of Probability} in 1986 \cite{MR866356} results under weak conditions on the first absolute moment of $\sum_{i=1}^n X_i$. As a consequence of one of their results, one can show that for $\phi$-mixing processes the central limit theorem holds if and only if
\[
  \limsup_{n\rightarrow \infty} \frac{\sigma_n}{ \E|\sum_{i=1}^n X_i| } \leq \sqrt{2/\pi}.
\]

These results from 1986 share a common feature, as they do not require any assumptions on the rate of convergence of the mixing coefficients. Besides motivating further refinements, they raise the question of to what extent a mixing condition is actually necessary for a central limit theorem.
In dynamical systems, this is a natural question, since many classical examples, such as an irrational rotation, are neither mixing nor weakly mixing in terms of ergodic theory. For example, Conze asked in the 1970s
whether there exists a probability-preserving dynamical system
$T: X \to X$ with zero entropy and a function $f : X \to R$ such that the stochastic process $(f \circ T^n)_{n\geq 0}$ satisfies the central limit theorem. This question was answered affirmatively by the Burton-Denker Theorem, published as ``{On the central limit theorem for dynamical systems}'' in \emph{Transactions of the American Mathematical Society}, in 1987 \cite{MR891642}, which tells us that such a function always exists for any aperiodic dynamical system.
With this result at hand, it is now natural to ask what are the  possible limit distributions of
$(f(X_n))_{n\geq 0}$, for a stationary prozess $(X_n)_{n\geq 0}$, not necessarily generated by a dynamical system, and some real valued function $f$. To be more precise, we say that $(f(X_n))_{n\geq 0}$ satisfies a distributional limit theorem if there are sequences $(a_n)_{n\geq 0}, (b_n)_{n\geq 0}$ and a random variable $Y$ such that
\[
  \frac{\sum_{i=1}^n f(X_i) - a_n }{b_n}  \stackrel{\mathcal{D}}{\longrightarrow} Y.
\]
In this context, the Burton-Denker theorem can be seen as a starting point for results in this direction. For example, Lacey (1991, \cite{MR1066446}) showed that any fractional Brownian motion can be obtained in the limit, and Liardet and Volny (1997, \cite{MR1459847}) showed that the limit distribution for all functions in a $G_\delta$ set are 
dense in the set of all probability measures, and that for a  rotation of the circle with a diophantine condition, this phenomenon also holds with respect to differentiable functions. On the other hand, Volny (1999, \cite{MR1624218}) and Lesigne (2000, \cite{MR1641132}) improved the  Burton-Denker theorem by showing that one may choose $f$ such that $(f \circ T^n)_{n\geq 0}$ satisfies an invariance principle and an almost sure CLT, respectively. 

The collaboration with Bob Burton (1948--2018), who had earlier obtained his Ph.D.\ with Don Ornstein at Stanford,  began when Burton spent the academic year 1983-84 as a Humboldt Research Fellow in G\"ottingen. This stay marked the beginning of a lifelong friendship with Manfred Denker and other members of Denker's scientific network. Until the mid-2000s, Burton spent several longer stays in Europe, either in G\"ottingen, or in Delft with Mike Keane, or in Groningen and later in Bochum with Herold Dehling. In return, Bob Burton's home base at Oregon State University in Corvallis became a favorite place for visitors from Europe.

\subsection{Limit theorems for  processes with strong mixing properties}
However, if the underlying process or dynamical system exhibits good mixing properties, one may seek stronger distributional limit theorems, such as local limit theorems, weak invariance principles, and almost sure versions thereof.
 For example, the almost sure invariance principle (ASIP) is an almost sure version of the weak invariance principle and is defined as follows. We say that the process $(X_n)$ satisfies the ASIP with rate $(r_n)_{n\geq 0}$ if, after eventually enlarging the underlying probability space, there is a Wiener prozess $(B(t))_{t\geq 0}$ and $r_n = o(\sqrt{n})$ such that
\[
\sum_{i=1}^n X_i =  B(n) + o(r_n)
\]
almost surely. Note that this is a very strong distributional limit theorem which has far-reaching consequences as almost every result known for the Brownian motion carries over to the partial sums of the $(X_i)$.

In the context of dynamical systems, the article ``{Approximation by {B}rownian motion for {G}ibbs measures and flows under a function}" by  Denker and  Philipp  (\emph{Ergodic Theory and Dynamical Systems} 1984, \cite{MR779712}) on flows built under a function is widely regarded as a breakthrough as their proof revealed how to employ $\psi$-mixing in order to establish the ASIP. Another highly influential contribution of Denker was his 1989 survey article ``{The central limit theorem for dynamical systems}"\cite{MR1102700}, which  established him as one of the leading experts in the field.

In this period, Denker also started working successfully with Jon Aaronson on non-integrable processes by importing ideas from infinite ergodic theory to this purely probabilistic setting. For example, they looked in ``{Lower bounds for partial sums of certain positive stationary processes}" (\cite{MR1035234}) into
 positive, non-integrable and $\psi$-mixing processes and were able to prove the other law of the iterated logarithm in this setting, meaning that
\[ \liminf_{n \to \infty} \frac{\sum_{i=1}^n X_i}{b(n/\log \log n) \log \log n }  \geq 1\]
almost surely with respect to a suitable normalising sequence $(b_n)_{n\geq 0}$, and that, if $b(n/\log \log n) \sim b(n)/\log \log n $, then the lims inferior is equal to one. 

A further, highly influential paper of these two authors, ``{Local limit theorems for partial sums of stationary sequences generated by {G}ibbs-{M}arkov maps}" from 2001 (\emph{Stochastics and Dynamics}, \cite{MR2048211}) uses a very careful adaption of Nagaev's method in order to show a distributional local limit theorem for Lipschitz observables in the domain of attraction of a stable law. In there, they considered a dynamical system $(X,T,m)$ with  the Gibbs-Markov property, which essentially means that $T$ admits a countable Markov partition, that $m$ is an invariant probability such that $\log \frac{dm}{dm\circ T}$ is Lipschitz continuous and that the associated transfer operator has a spectral gap. For $\phi: X \to \R$ Lipschitz continuous and aperiodic in the domain of attraction of a stable law of order $0< \rho < 2$, they showed that there exist $(a_n)_{n\geq 0}$, $(b_n)_{n\geq 0}$ such that $((\sum_{i=1}^n   f\circ T^i  -  a_n)/b_n)_{n\geq 0}$ converges in distribution to a stable law of order $\rho$. Moreover, the local version of this distributional limit holds: Assume that $(k_n - a_n)/b_n \to \kappa$ and that $a < b$. Then
\[
\lim_{n \to \infty} b_n m \left(\left\{ x : a +  k_n  \leq   { \textstyle \sum_{i=1}^n }  f\circ T^i(x)  \leq b +  k_n  \right\}  \right)
=  (b-a) f_\rho(\kappa),
\]
where $f_\rho$ is the density of the $\rho$-stable distribution. Besides the obvious interest in this local limit theorem, the paper became very influential due to its applicability to Young towers and its method of proof.

But Denker did not stop looking for new and interesting open questions. For example, one could ask for almost sure versions of local limit theorems, which is a new field of research with relatively few examples and publications. However, one of those is the paper by Susanne Koch and Denker from 2002, entitled ``{Almost sure local limit theorems}" (\emph{Statistica  Neerlandica}, \cite{MR1916315}), which establishes this kind of local limit theorem for the symmetric random walk. 

On the other hand, instead of considering the partial sum of a single function along an orbit, one may consider a sequence  $(f_n)_{n \geq 0}$ of functions and the sum $\sum_{k=0}^{n-1} f_k \circ T^k$, which is a canonical object for random and sequential dynamical systems. For a fixed dynamical system, this is related to the very active fields  of research on shrinking targets and extreme value theory.  In ``{A Poisson limit theorem for toral automorphisms}" (\emph{Illinois Journal of Mathematics} 2004,  \cite{MR2048211}), Manfred Denker, Mikhail Gordin and Anastasya Sharova took a different approach and considered an abstract sequence of sets $(G_n)_{n\geq 0}$ with shrinking diameter, 
whose set of accumulation points contains no periodic point of the underlying toral automorphism. Then, for $(s_n)_{n \geq 0}$ and $\lambda> 0$ with $\lim_{n \to \infty} s_n \mathrm{Leb}(G_n) = \lambda$, under some additional regularity condition on $(G_n)_{n\geq 0}$, 
\[
\sum_{i=0}^{s_n-1} \mathbf{1}_{G_n} \circ T_i 
\]
converges to a Poisson distribution with parameter $\lambda$.  To give an explicit example, one may consider a sequence of balls $(B_{r_n}(x_n))_{n\geq 0}$ such that $r_n \to 0$ and the accumulation points of $(x_n)_{n\geq 0}$ do not contain a periodic point. Here, it is worth noting that the condition of non-accumulation at periodic points is necessary, as shown by Hirata (1993, \cite{MR1245828}), and that the paper relies on the Chen-Stein method and Gordin's homoclinic Laplace operators, which makes the proof very interesting. Probably due to these characteristics, the paper was well received by the community.

One of the most recent papers by Denker on probabilistic limit laws for dynamical systems makes use of one of the most prominent objects in hyperbolic dynamics: periodic points. In the paper     
``{Fluctuations of ergodic sums of  periodic orbits under specification}" by Manfred Denker, Samuel Senti and Xuan Zhang (\emph{Discrete and Continuous Dynamical Systems} 2020, \cite{MR4112026}), the authors show how uniform measures on periodic orbits can be used to obtain a Lindeberg-type central limit theorem with respect to the measure of maximal entropy. This appears to be a new observation and was subsequently further developed by Thompson and Wang.

\subsection{Perspective}

Denker’s work in probability theory is characterized by the analysis of stochastic processes with strong dependence, often arising from dynamical systems. A central theme is the extension of classical limit theorems to dependent settings, together with the identification of minimal conditions under which such results hold.
In this context, Denker’s contributions range from central limit theorems under weak assumptions to refined limit laws for strongly mixing systems. At the same time, results such as the Burton--Denker theorem show that probabilistic limit behavior can occur even in the absence of mixing, thereby significantly broadening the scope of limit theory for dependent processes.

\section{Fractal Sets}
\label{sec:fractals}

Fractal geometry and dimension theory are central and sustained areas of Manfred Denker’s research, particularly from the late 1980s onward, with continuing developments in later work on iterated function systems and probabilistic boundaries. A characteristic feature of this work is an extensive and deep use of ergodic theory and the thermodynamic formalism to relate geometric properties of fractal sets to invariant measures and its statistical structure.

A recurring theme is the translation of geometric information, such as Hausdorff dimension or boundary behavior, into measure-theoretic terms via invariant or conformal measures. This perspective allows tools from probability and dynamical systems to be applied to problems in fractal geometry.

 Denker has been one of the pioneers in developing the theory of conformal measures, which has shaped to today the broader program of expressing geometric invariants of dynamical systems in terms of thermodynamic quantities.

\subsection{Conformal measures and dimension theory for rational maps}

A major focus of Denker’s work on fractal sets during the period 1986--1991,
largely in collaboration with Mariusz Urba\'nski, was the development of a
systematic framework to treat questions concerning the fractal dimensions (Hausdorff, Box counting, Packing, etc.) for rational maps and related classes of 
dynamical systems.  Two papers from the early 1990's play a key role in this program.

In the 1991 {\em Nonlinearity} paper titled ``{On Sullivan’s conformal measures for rational maps of the Riemann
sphere}"  \cite{MR1107011}, Denker and Urba\'nski establish a precise
connection between certain thermodynamic quantities and the geometric dimension.  They show that for a rational map $T$, the supremum of the Hausdorff dimensions of all ergodic, $T$-invariant probability measures with positive entropy coincides with the minimal zero of the pressure function of the geometric potential associated to $T$. The supremum is nowadays called the dynamical dimension of the system. More precisely, if $P(t)$ denotes the pressure function, then the Hausdorff dimension is characterized by the unique value $t$ for which $P(t)=0$. 
While similar results were established earlier by Bowen and Ruelle \cite{MR556580,MR684247} for hyperbolic maps,
the crucial improvement was that Denker and Urba\'nski's result holds for all expansive rational maps. The thermodynamic dimension formula for parabolic maps was later generalized by Urba\'nski and Christian Wolf to certain diffeomorphisms in $\mathbb{R}^2$ with parabolic behavior, so-called parabolic horseshoes, \cite{UW2}.
Denker's pioneering result firmly connects Hausdorff
dimension to the framework of thermodynamic formalism and clarify
the role of conformal measures which were introduced earlier by  Sullivan, in describing the finer structure of Julia sets and more general invariant set.
Closely related ideas appear in 
``{Ergodic theory for Markov fibred systems and parabolic rational maps}"
(1993) \cite{MR1107025}, where Aaronson, Denker and Urba\'nski analyze dynamical systems with
parabolic behavior.  Among other results, they prove that the Hausdorff dimension
of a parabolic Julia set in the Riemann sphere is strictly less than two, a fact
that highlights the subtle geometric consequences of non‑hyperbolic dynamics.
Taken together, these works provide a fairly complete picture of how different types
of dynamical behavior, hyperbolic, parabolic, affect the dimension and
measure-theoretic structure of the associated fractal sets.

These works exemplify a central goal of Denker during this time period: to unify
various notions of dimension and to relate them quantitatively to entropy and
pressure.  A broader perspective on this program and its connections to earlier
work of Patterson and Sullivan can be found in the later survey by Stratmann and
Denker {``Patterson measure: Classics, variations and applications''}
(2012) \cite{MR3060460}.

\subsection{Existence of conformal measures for iterated function systems}

Denker's dimension--theoretic results for rational maps are closely related to an
earlier foundational contribution,
``{On the existence of conformal measures}" (1991) \cite{MR1014246}, in which
Denker and Urba\'nski establish the existence of conformal probability measures
for general countable iterated function systems and families of potentials.
Motivated in part by Patterson’s approach to Fuchsian groups \cite{MR450547},
the paper shows how
geometric expansion properties can be encoded by measures that transform
conformally under the action of the system.

While conformal measures had appeared earlier in work of Sullivan (who introduced the term ``conformal density" in \cite{MR556586,MR634434}), Rohlin, and Keane (in the context of $g$-measures),
the contribution of Denker and Urba\'nski lies in developing a general and flexible framework that connects these constructions directly to geometric and dynamical properties. In particular, conformal measures provide a tool for encoding metric expansion in measure-theoretic terms.

These ideas were revisited and significantly generalized more than two decades
later in
``{Conformal families of measures for general iterated function systems}"
(2015) \cite{MR3330340}, joint with Yuri.  In this work, the authors extend the
theory to systems of homeomorphisms acting on compact metric spaces without any
assumed Markov structure, encompassing graph‑directed Markov systems and inverse
branches of more general dynamical systems.  The introduction of conformal
\emph{families} of measures allows finer control of boundary behavior and
describes how conformal measures interact with the geometry of the underlying
space.

Seen together, the 1991 and 2015 papers highlight the lasting importance of
conformal measures in Denker’s work on fractal geometry: they provide a
conceptual bridge between geometric expansion, thermodynamic formalism, and
measure--theoretic structure.

\subsection{Fractal boundaries and the Sierpiński gasket }

A closely related, but different, aspect of Denker’s work on fractal sets appears
in his study of probabilistic and potential--theoretic representations of
classical self‑similar fractals.  In their paper
``Sierpi\'nski gasket as a Martin boundary I. Martin kernels''
(\emph{Potential Analysis} (2001), \cite{MR1822915}), Denker and Sato show that the
Sierpinski gasket can be realized as the Martin boundary of an explicitly
constructed Markov chain.

This work was the first in a sequence of papers, written in collaboration
with Sato and later with Koch between 1999 and 2001 (\cite{MR1763899,MR1739300,MR1822915,MR1912376}), devoted to the realization of
the Sierpi\'nski gasket and related fractals as boundaries of stochastic processes.
The approach combines ideas from potential theory, Markov chains, and fractal
geometry, and provides a probabilistic interpretation of harmonic functions on
self--similar sets.

By representing the gasket as a Martin boundary, Denker’s work connects the
analytic structure of the fractal to probabilistic behavior at infinity, adding a
new layer to the study of classical fractals.  This line of research
illustrates how fractal geometry, probability, and ergodic theory interact in a
natural and conceptually unified way.

\subsection{Perspective}

Denker’s contributions to fractal geometry form part of a broader program linking dynamical, probabilistic, and geometric viewpoints. Central to this program is the use of invariant and conformal measures to encode geometric complexity and to relate dimension, entropy, and pressure.

From rational maps and parabolic dynamics to general iterated function systems and probabilistic boundaries, these ideas have provided a unifying framework for understanding fractal structures in dynamical systems. They continue to influence current research on thermodynamic formalism and dimension theory.

\section{Dynamical Systems}
\label{sec:dynamical_systems}

In this section we survey some of the main contributions of Manfred Denker to the theory of dynamical systems. The list of results had to be reduced to a conservative selection, which was suggested by Denker himself. The methods applied are primarily from measurable dynamics and ergodic theory, and often involve the analysis of the pressure function associated with the underlying dynamical system. Denker's main collaborators were Walter Philipp,
Mariusz Urba\'nski (a postdoctoral assistant of Denker in G\"ottingen during the late 80's and early 90's), Jon Aaronson and Mikhail Gordin. From 1984 on, the Sonderforschungsbereich (SFB) 170 in Geometry and Analysis provided Denker with funds for several doctoral and postdoctoral research positions, as well as for support of short term visitors. 
During the SFB years  Denker and Samuel Patterson were administrating the project on Dynamical Systems and Analytic Number Theory which included a popular weekly research seminar. 
 Throughout his career Denker enjoyed collaborating across all mathematical maturity levels---doctoral students, postdocs and established mathematicians---many times bringing all of those together to a single project. To select a few examples: from Denker's early days in G\"ottingen 
 with several collaborations with Dehling and Philipp, over the many papers with Urba\'nski and Aaronson, to the 2017 paper with Zhang and Senti  \cite{MR3734148} and the Ergodic Theory and Dynamical Systems paper with Stadlbauer and Kifer \cite{MR2380303}. 
One of Denker's particularly productive periods was between the late 1980's and early 1990's, when Mariusz Urba\'nski was an SFB postdoc and Jon Aaronson a regular visitor. Several papers from this period fall into the realm of Dynamical Systems and will be discussed below.

\subsection{ Approximation by Brownian motion for flows}
While most of the publications of Manfred Denker with Walter Philipp had statistical and probabilistic objectives, and for that part are all with Herold Dehling, there is one paper with a dynamical focus, namely    
 ``{Approximation by Brownian motion for Gibbs measures and flows under a function}" published 1984 in {\em Ergodic Theory and Dynamical Systems} \cite{MR779712}. In the paper the authors employ methods from probability theory, namely the Wiener Process, aka Brownian motion, to approximate (geodesic) flows defined by a roof function over a topological Markov chain equipped with a Gibbs measure. As consequence of the main approximation theorem the authors derive a 
 central limit theorem (CLT), a weak invariance principle and the law of the iterated logarithm. This result extends the CLT of Ratner \cite{zbMATH03424864} for geodesic flows on negatively curved manifolds, obtained by coding the geodesic flow as a flow over a shift space, and shows that several fundamental statistical properties hold. It is worth noting that several forms of mixing, among those $\phi$- and $\psi$-mixing, as displayed in the probability section of this article, play a role in the main theorem.

Consider the flow $S_t$ under a H\"older continuous function $l: \Sigma \rightarrow \mathbb{R}_{+}$ with base transformation $(\Sigma, T, \mu)$, where $\Sigma$ is a topological Markov Chain with shift $T$ and $\mu$ a Gibbs measure. Let  $X=\{(x,s): x \in \Sigma ,\ 0 \leq s < l(x)  \}$ be the space on which $S_t$ is defined,  
equipped with probability measure $\nu$ given by $ d\nu=(\int l \ d \nu)^{-1}\ ds \times d\mu$.
If $f$ is a measurable real valued function on $X$ with $Ef =0$, $E|f|^{2+\delta}<\infty$ for 
some $0<\delta \leq 1$ and so that  $E|f-E[ f|(\gamma)_{-n}^{n}]|^{2+\delta} \ll n^{-(2+7/\delta)(2+\delta)}$ where 
$(\gamma)_{-n}^{n}$ is the $\sigma$-field generated by the sets of the form $\{(x,s):\ x \in A, \ 0 \leq s <l(x) \}$ and  
$A$ is a central cylinder of length $2n+1$. Under these assumptions on $f$ Theorem 1.1 in \cite{MR779712} 
states that 
\[ \sigma^2=\lim_{t \to \infty} \frac{1}{t} E\left( \int_0^t f(S_{\tau}(w)) \ d\tau\right)^2\] 
exists and if $\sigma^2>0$ there exists a standard Brownian motion $\{B(t,x): t \geq 0,\ x \in \Sigma\}$ on $(\Sigma, \mu)$ so that 
\[ \sup_{0 \leq u <l(x)} \left| \int_0^t f(S_{\tau}(x,u)) d\tau -B(\sigma^2 t, x) \right| \ll t^{1/2 - \lambda} \] 
for $\mu$ almost all $x$ and $0<\lambda<\delta/588$.

Denker and Philipp draw three important conclusions, besides the theorem invoking only standard results.  
A {\em central limit theorem:} For all $z \in \mathbb{R}$
\[ \lim_{t \to \infty} \nu\left \{w:\ \frac{1}{\sigma \sqrt{t}}\int^t_0 f(S_{\tau}w)\ d\tau \leq z \right \} =\frac{1}{\sqrt{2\pi}} 
\int^z_{-\infty} \exp\left(-\frac{1}{2}u^2\right)\ du.\]
{\em The law of the iterated logarithm:} For $\nu$-almost all $w$: 
\[  \lim_{t \to \infty} \frac{1}{\sqrt{2\sigma^2 t\log\log t}} \int^t_0 f(S_{\tau}w)\ d\tau=1.\] 
and {\em The weak invariance principle:} Let
\[ \zeta_n(t)=\frac{1}{\sigma \sqrt{n}}  \int^{nt}_0 f \circ S_{\tau} \ d\tau \quad (t \in [0,1], n \in \mathbb{N}), \]
then $\{\zeta_n\}_{n \in \mathbb{N}}$ converges weakly in $C([0,1])$ to the standard Wiener measure. \\

The paper received a lot of interest, often times as motivational and foundational.  
It has been cited, for example, by Daniel Rudolph in his 1991 publication  ``{Asymptotically Brownian skew products give non-loosely Bernoulli K-automorphisms}" \cite{zbMATH04069892} 
and seventeen years later by Curtis T.\ McMullen in ``{ Thermodynamics, dimension and the Weil-Petersson metric}" 
\cite{zbMATH05313092}. 
Generally, the central limit theorem and its variants and applications to dynamical systems is one of the icons in Denker's research and naturally appears in several parts of this article. To this end we mention the 1987 Transactions paper ``{On the central limit theorem for dynamical systems}" with Robert Burton \cite{MR891642}. 
In fact, Denker typically proved the central limit theorem in his probability lectures for undergraduates and it was not uncommon that students attempted to prove it as part of their oral pre-diploma examinations with Manfred Denker.

The collaboration with Walter Philipp (1936--2006) was initiated by the joint supervision of the Ph.D.\ research of Denker's student Herold Dehling, who spent the academic year 1979-80 as a visiting graduate student in Philipp's lab at the University of Illinois at Urbana-Champaign. Subsequently, Walter Philipp spent several summers in G\"ottingen, 
where most of the joint papers  with Dehling and Denker were written. Their collaboration 
resulted in a life-long friendship that ended only with Philipp's untimely death in 2006.

\subsection{Hopf decomposition}
The Hopf decomposition theorem is one of the fundamental results in the study of dynamical systems, as it allows one to decompose the underlying phase space into a conservative and a totally dissipative part. Furthermore, as shown by Hopf, this decomposition is trivial for the geodesic flow on a manifold of constant negative curvature. Moreover, in this setting, he observed that conservativity implies ergodicity. In short, the flow is either totally dissipative, or conservative and ergodic, which is one of the earliest results on geodesic flows on manifolds.

However, for general non-singular transformations, this dichotomy fails for several reasons. First, the decomposition into conservative and dissipative parts may be non-trivial. Moreover, even when the decomposition is trivial, all combinations can occur: there exist conservative but non-ergodic transformations, as well as ergodic, totally dissipative and non-invertible transformations.

In particular, an analysis of Markov fibred systems should start with a discussion of this dichotomy, as in ``{Ergodic theory for Markov fibred systems and parabolic rational maps}" by Aaronson, Denker and Urba\'nski from 1993 (\emph{Transactions of the American Mathematical Society}, \cite{MR1107025}). There, the authors prove, among many other important results, the precise analogue of Hopf's dichotomy under suitable weak distortion assumptions. That is, a topologically transitive Markov fibred system is either totally dissipative or conservative and ergodic. We will address this paper in more detail below. 

Furthermore, the analogue also holds for random Markov fibred systems, as shown by Denker,  Kifer and Stadlbauer  in ``{Conservativity for random Markov fibred systems}", \emph{Ergodic Theory and Dynamical Systems} from 2007 \cite{MR2380303}. However, it should be noted that in the random case the proof is considerably more involved, as the authors are required to work with a generalized class of random Markov fibred systems in order to control first return maps to cylinders. This extension is of independent interest, as it allows the application of many standard techniques from dynamical systems. Moreover, this dichotomy proves to be an essential ingredient in the proof of the random Perron-Frobenius-Ruelle theorem for countably many states in \cite{MR2410952}. 

\subsection{Hausdorff and conformal measures on Julia sets with a rationally indifferent periodic point.}
In Denker and Urba\'nski's 1991 \emph{Journal of the London Mathematical Society} paper \cite{MR1099090} bearing this title, they consider holomorphic endmorphisms  $T: \overline{ \mathbb{C} }\to \overline{\mathbb{C}}$ (of degree $d\geq 2$) on the Riemann
sphere $\overline{ \mathbb{C}}$ and their Julia set $J(T)$. For $t \geq 0$ a probability measure $m$ on $J(T)$ is called $t$-conformal for $T: \overline{\mathbb{C}} \rightarrow \overline{\mathbb{C}}$, if 
\[ m(T(A))=\int_A |T'|^t\ dm\] 
for every Borel set $A \subset J(T)$ satisfying that $T|_A$ is injective. Sullivan showed that for {\em expanding} (hyperbolic) $T$, there is a nonatomic and unique $\delta$-conformal measure where $\delta=HD(J(T))$ is the Hausdorff dimension 
of $J(T)$. 

Denker and Urba\'nski give a similar characterization of $HD(J(T))$ for {\em expansive} $T:J(T)\rightarrow J(T)$. In fact, 
they show that $T|_{J(T)}$ is expansive if and only if $J(T)$ contains no
critical points of T (Theorem 4). The authors conclude that any expanding mapping $T: J(T) \rightarrow J(T)$ is expansive and that an expansive $T$ is not expanding, if and only if $J(T)$ contains an indifferent rational periodic point. 
They further show that for a rational map with a rationally indifferent periodic point (and no critical points in the Julia set), the Hausdorff dimension $HD(J(T))$ of $J(T)$ equals the smallest zero of the pressure function $t  \to P(T,-t \log |T'| )$ and that $t=HD(J(T))$ is the smallest exponent such that a $t$-conformal measure exists (Theorem 15). 
In addition they obtain the formula 
\[  HD(J(T))=\sup_{\mu } \inf_{Y \subset J(T)} \{ HD(Y):\ \mu(Y)=1\},    \]
where  $\mu$ ranges over the $T$-invariant, ergodic, Borel probability measures with positive entropy. 
From the many citations this paper received we like to mention Ian Melbourne and Dalia Terhesiu with their 2012 {\em Inventiones Mathematicae} paper \cite{MR2929083}.

\subsection{Hopf decomposition for Julia sets} 
Here we present a broad account of the paper ``{Ergodic theory for Markov fibred systems and parabolic rational maps}"   of Jon Aaronson, Manfred Denker and Mariusz Urba\'nski \cite{MR1107025}.  
This paper acts as a bridge between thermodynamic formalism and the geometric measure theory of non-uniformly hyperbolic rational maps, extending tools previously restricted to uniformly expanding maps to parabolic cases.
Previous work of Denker and Urba\'nski \cite{MR1129999} 
shows that for a parabolic rational map $T: \overline{\mathbb{C}}\rightarrow \overline{\mathbb{C}}$ with $h$-conformal measure $m$, there is a topological Markov partition with respect to which $(J(T), m, T)$ is a fibred system. As before  $J(T)$ denotes the Julia set of $T$. 
The key mathematical tool developed in this work are {\em (Markov) fibred systems}. For the reader's convenience, we recall their definition. 

Let $T$ be nonsingular on the nonatomic probability space $(X, \mathcal{F}, m)$. 
A (Markov) fibred system on $X$ is a pair $(T, \mathcal{R})$  where $ \mathcal{R}\subset  \mathcal{F} $ is a countable partition of $X$ that generates $\mathcal{F}$ and has the {\em Markov property}, that is  
\[ TB = \bigcup_{B' \in \mathcal{R}, \ B' \cap TB \neq \emptyset}   B' \quad \text{ for all } B \in \mathcal{R}.\]
In addition, $T$ is required to be {\em locally invertible}, i.e. for every $B \in \mathcal{R}$, $m(B)>0$ and $T|_B$ is nonsingular and invertible. 

In the first part of the paper the authors derive general properties of Markov fibred systems, in particular a Hopf dichotomy.   Namely, if $(T,\mathcal{R})$ is an irreducible fibred system with the Schweiger property, then $T$ is either conservative or totally dissipative and if $T$ is conservative, then $T$ is ergodic. Moreover, if $\inf \{m(B):\ B \in T\mathcal{R} \} >0$, then $T$ is conservative and ergodic (Theorems 2.5. and 2.8). It is further shown that the theory applies to  Markov shifts and certain piecewise linear maps of the unit interval. 

In section 3 of the paper, the authors show that for an irreducible fibred system $(T,\mathcal{R})$ that has the Schweiger property with respect to $\mathcal{R}(C,T)$ (for a detailed definition of these set systems see the paper) and conservative $T$ there is a $\sigma$-finite, $T$-invariant measure $\mu \sim m$ whose logarithmic  Radon-Nikodym derivative is bounded a.e.\ on sets in $\mathcal{R}(C,T)$. Furthermore any set  $\mathcal{R}(C,T)$ is a  Darling-Kac set for $T$ whose return time process is continued fraction mixing. It is also noted that for conservative  fibred systems induced systems are canonically fibred. This section also shows that $T$ is exact whenever $(T,\mathcal{R})$ is an aperiodic fibred system with the Schweiger property and conservative $T$.  
  
In section 6 of the paper mixing properties of fibred systems with finite invariant measures are studied. The authors prove if $(T,R)$ is aperiodic and fulfills some assumptions not detailed out here, there exists a $T$-invariant probability measure $q \sim m$, so that $(T,R)$ is continued fraction mixing. Not all that surprising, the paper contains a CLT for certain fibred systems (Theorem 7.2). The remainder of the paper deals with parabolic rational maps $T : \overline{\mathbb{C}} \rightarrow  \overline{\mathbb{C}}$. Here the existence and uniqueness of an h-conformal measure for $T$ is shown and that such a measure is non-atomic (Theorem 8.7). Under the same assumptions on $T$ it is also shown that $HD(J(T))<2$ (Theorem 8.8).   This result was later re-proved by Lukas Geyer who established that parabolic Julia sets are porous which implies $HD(J(T))<2$ \cite{MR1717570}.

Finally, the theory is applied to parabolic rational maps of the Riemann sphere. More specifically, the situation where  $T$ is a non-singular transformation on the probability space $(J(T), \mathcal{F}, m)$ and $m$ is the h-conformal measure on the Borel sigma algebra $\mathcal{F}$ on $J(T)$. The authors show first that there is a partition $\mathcal{R}$ of $T$ such that $(T,\mathcal{R})$ is a finite, aperiodic, parabolic fibred system with aperiodic jump transformation. Then the established results imply that there is a finite invariant measure $q \sim m$, that has bounded and uniformly Lipschitz Radon-Nikodym derivative (Theorem 9.6). There are more theorems that essentially display and refine the general results for fibred systems to parabolic rational maps. Among them are two versions of central limit theorems (Theorems 9.12 and 9.13).   

This is one of Denker's most cited papers, and includes citation from Matthias Meiwes 2025 \cite{MR4966568},  Yuri Lima and Omri Sarig 2019 \cite{MR3880208},  
Ian Melbourne and Dalia Terhesiu 2012 \cite{MR2929083}
and from Artur Avila together with Misha  Lyubich 2008 \cite{MR2373353}.

\subsection{Gibbs measures for fibred systems} 
In Denker's joint article with Mikhail Gordin  in {\em Advances in Mathematics} (1999) \cite{MR1736956},
they consider Gibbs families of conditional measures for fibred systems by looking at relative 
transfer operators. More precisely assume $T: Y \rightarrow Y$ and $S: X\rightarrow X$ are continuous maps of 
Polish spaces, the latter being a factor of the first via a continuous factor map $\pi: Y \rightarrow X$. 
If $T$ is fibrewise expanding and exact along fibres, i.e., maps fibres onto fibres,  and if $\phi$ is a H\"older continuous function the existence of the system of conditional measures $\{\mu_x\}_{x \in X}$, called a family of Gibbs measures, on the fibers $Y_x=\pi^{-1}(x)$ is shown. The fibre-maps are assumed to have universally bounded preimages $\# T^{-1}_x(y) \leq M$ for all $x \in X$ and $y \in Y_x$. Denker's and Gordin's assumptions on the base transformation differs from the ones made in previous works as they consider continuous, but not necessarily invertible, onto maps. 
Noninvertibility presents difficulties, making the results differ from the classical Ruelle theorem. 
While the existence and uniqueness of Gibbs measures still holds, the part on invariant measures does not generalize, but sufficient conditions for their existence are given. The construction of Gibbs families for fibred systems uses Birkhoff 's theory and not symbolic representations as one may expect.

In some more detail, a system $\{\mu_x: \ x \in X\}$ of conditional probabilities for the fibred space at hand is called a {\em Gibbs family} for a measurable $\varphi: Y\rightarrow \mathbb{R}$, if there exists a positive, measurable function $A: X \rightarrow \mathbb{R}$ so that for all $x \in X$ the Jacobian of $\mu_x$ with respect to $T$ is given by 
\[ \frac{d\mu_{S(x)}\circ T_x}{d\mu_{x}}=A(x)\exp[-\varphi], \quad \mu_x \ a.e.\]
Following this definition, Denker and Gordin are able to show several theorems, the first one (Theorem 2.6) addresses the {\em existence} of Gibbs measures: { Given a compact fibred system with a metric $d$ on $Y$, that is fibrewise expanding and topologically exact along fibers, then for every H\"older continuous function $\varphi: Y \to \mathbb{R}$ there exists a unique Gibbs family $\{\mu_x: x \in X\}$ for $\phi$, so that $\text{supp}(\mu_x) = Y_x$ for all $x \in X$. Given $\varphi$, the scaling function $A$ is also unique.}  

Next, the authors add a continuity result (Theorem 2.7), where the map 
$i: Y \to \{(x,y) \in X\times Y: \  S(x)=\pi(y)\}$ given by $i(y):=(\pi(y), T(y))$ is employed: 
{ If $\{\mu_x: x \in X\}$ and $(\varphi,A)$ are as above, $S$ and $\pi$ are open maps and $i$ is a local homeomorphism, then $\{\mu_x: x \in X\}$ and $A:X\rightarrow \mathbb{R}$ are continuous.} 

The third fundamental result (Theorem 2.10) shows H\"older continuity of the maps $x \mapsto \int_Y f(y)\ d\mu_x(y)$ 
and $x \mapsto A(x)$ under the assumptions that $T$ is open, fulfills a H\"older condition, $\pi$ a contraction and 
$S$ expansive, Lipschits continuous and locally non-contracting. 

As an application,  Denker and Gordin first show that if $S$ is expansive, every H\"older continuous 
function $\phi$ admits a Gibbs measure $\mu$ on $Y$, so that its family of conditional measures is the unique 
Gibbs family for $\phi$ supported on the leaves of the foliation (Theorem 3.1). If additionally the 
system $(Y,T)$ is expanding  there exists a measure $\nu \ll \mu$, so that its family of conditional measures 
$\{\nu_x\}$ is of the form $\nu_x=h_x \mu_x$ almost everywhere with a continuous $h_{\pi}: Y \rightarrow \mathbb{R}$ (Theorem 3.2). With these preparations Denker and Gordin apply their results to polynomial endomorphisms $\widehat T: \mathbb{C}^2 \to \mathbb{C}^2$ and its restriction $T$ to its Julia set $J(\widehat T)$ (Example 3.3). In this case the fibration is defined by the projection of $\widehat T$ to the first factor. Then (Theorem 3.4) the existence of a unique equilibrium measure $\mu$ to any given H\"older continuous function $\phi$ is shown. Moreover it is shown that $\mu$ admits a unique H\"older continuous disintegration that is absolutely continuous with respect to the Gibbs family for $\phi$.

\subsection{Lindeberg theorem for Gibbs-Markov dynamics}
This so-entitled article appeared in {\em Nonlinearity} (2017), with Samuel Senti and Xuan Zhang, \cite{MR3734148}.
This paper presents (briefly) a first significant step into the study of the Central Limit Theorem (CLT) in dynamical systems for arrays of functions. This is the analog of an array known in probability and the Lindeberg condition. In the paper it is shown that such arrays can be handled in Gibbs-Markov dynamical systems to derive the CLT. Particularly, it is shown that the space of all $L^2$ functions which satisfy the CLT is much larger than the class of H\"older continuous functions. The result is important for considering inference properties of Dynamical Systems.  On the other hand, it is still an open problem to determine the class of functions in $L^2$ satisfying the CLT. For those reasons, the paper opens up a wide area of further research.

For a {more detailed summary}, recall that Lindeberg's central limit theorem deals with arrays of independent random variables, i.e., families of random variables defined on different probability spaces. This paper formulates Lindeberg's central limit theorem for {\em dynamical} arrays, and proves such a theorem for arrays in dynamical systems, such as Gibbs-Markov maps, or systems as formulated by Aaronson and Denker 2001 in ``{Local limit theorems for partial sums of stationary sequences generated by Gibbs-Markov maps}" \cite{MR1840194}. Examples include certain countable state Markov chains and Markov maps of the unit interval and parabolic rational maps.    
Such theorems provide central limit theorems for Birkhoff sums $\sum^n_{i=1} f \circ T_i$ 
for certain functions which are not Lipschitz-, or H\"older-continuous. 

Theorem 4.1 (also see \cite{zhang}) presents a CLT for certain function sequences on Gibbs-Markov systems. 
Known examples of general settings to which this theorem applies, beyond the Gibbs-Markov systems, 
include maps of the interval endowed with the bounded variation norm, as well as Young towers endowed with the H\"older norm. Continued fraction transformations with Gauss measure also fulfill the CLT.

Section 5 contains the CLT for dynamical arrays.  To start,  
a dynamical array is a sequence $\{(F_{n,i}, \tau_{n,i}):\ i = 1, ..., k_n \}_{n \in \mathbb{N}}$
consisting of a family of real valued functions $F_{n,i}$ defined on a dynamical system
$(\Omega, T )$ and a family of initial times $\tau_{n,i} \in \mathbb{N}$, where $F_{n,i}$ is of the form 
\[  F_{n,i}=\sum^{l_{n,i}}_{j=1} f_{n,i,j}\circ T^{j-1}  \]
with $f_{n,i,j}: \Omega \to \mathbb{R}$, $l_{n,i} \in \mathbb{N}$ and $\tau_{n,i-1}+l_{n,i-1} \leq \tau_{n,i}$ 
for all $i=2,...,k_n$. 

Now for an array $\{(F_{n,i}, \tau_{n,i}):\ i = 1, ..., k_n \}_{n \in \mathbb{N}}$ over a Gibbs-Markov system 
$(\Omega, \mu, T,\alpha)$, with  $F_{n,i}$ as defined above, with centered $f_{n,i,j} \in L$.   
The function space $L$ is a subspace of $L^{\infty}$ equipped with a H\"older norm, subject to the partition given  
by the Gibbs-Markov system. 
Let $A_{n}:=F_{n,1}\circ T^{\tau_{n,1}}+...+F_{n,k_n}\circ T^{\tau_{n,k_n}}$.
If $\text{Var}(A_{n})>0$, a Lindeberg condition 
holds and in addition some finiteness conditions regarding the $\tau_{n,k_n}$ hold, the array satisfies the CLT 
(Theorem 5.3):   
\[  \frac{A_n}{\text{Var}(A_{n})} \Rightarrow   N(0,1).\]
This theorem has multiple applications in statistics and section 6 of the paper presents one: The Behrens-Fisher problem to determine whether two distributions are different, or not, in a statistical sense. To do this in the non-parametric setup, the authors show asymptotic normality of the Wilcoxon two sample rank statistics. 
More examples are outlined in Zhang's thesis \cite{zhang} ``{Studies on the weak convergence of partial sums in Gibbs-Markov dynamical systems}" (2015) supervised by  Denker.

\subsection{On specification and measure expansiveness.} 
This article appears in {\em Discrete and Continuous Dynamical Systems} (2017) with Welington Cordeiro and Xuan Zhang, \cite{MR3640582}.

While this contribution seems topically and methodically rather nontypical for Denker's work, it actually 
roots in his first interests as a young mathematician. Namely, in Denker's diploma thesis he presented the early work of R.\ Bowen on Axiom-A-diffeomorphisms with some improvements. The results were presented at a conference in Paimpol 1970 \cite{paimpol70}. In particular, as a novelty, the existence of Markov partitions for Axiom-A-homeomorphisms was shown, a notion which later was termed by Ruelle as {\em Smale spaces}.  One of the basic results in Bowen's work is the proof of Smale's spectral decomposition theorem (Theorem 2.5 in \cite{paimpol70}). In the work at hand, it is shown that specification is sufficient to obtain such a result (see Theorem D below).  This particular work presents an example of  
Denker working with two young mathematicians at the beginning at their career.  
The following theorems, among them the main Theorems A, B and C, display the interrelation of various 
concepts of shadowing and specification:  \\

\noindent In Theorem 3.2, $f$ has the local weak specification property if and only if $f$
has the shadowing property.\vspace*{2mm}

\noindent Theorem A: 1.\ $f$ has the limit shadowing property if and only if $f$ has
the local limit specification property. 2.\ If $f : X \to X$ is a Lipschitz map, then $f$ has the Lipschitz shadowing
property if and only if $f$ has the local Lipschitz specification property. 3.\ $f$ has the two-sided limit shadowing property if and only if $f$ has the local two-sided limit specification property.\vspace*{2mm}

\noindent Theorem B: If $f$ has the local specification property, then $f$ has the periodic shadowing property.\vspace*{2mm}


\noindent Theorem C: Let $f : X \to X$ be a transitive and strong measure expansive
homeomorphism. Then the following equivalences hold: 
$ f$ has the local specification property $\Leftrightarrow$
$ f$ has the local weak specification property $\Leftrightarrow$
$ f$ has the shadowing property $\Leftrightarrow$
$ f$ has the strong periodic shadowing property $\Leftrightarrow$
$ f$ has the periodic shadowing property $\Leftrightarrow$
$ f$ has the special shadowing property. \vspace{2mm}

\noindent The following Spectral Decomposition Theorem (Theorem D) develops ideas already contained in Denker's diploma thesis.  \vspace{2mm}
 
\noindent Theorem D: Let $f$ be a strong measure expansive homeomorphism with shadowing, then: 
 
\noindent 1. $ \Omega(f) = \bigcup^l_{i=1} B_i$ where $B_i$ is closed, invariant and $f|_{B_i}$ is transitive. 

\noindent 2. For each $0 \leq k \leq l$ $B_k = \bigcup^{\alpha_k -1}_{i=0} C_i$, where $C_i \cap C_j = \emptyset$ if $i \neq j$, 
and for each $0 \leq i < a_k$, $C_i$ is closed, $f(C_i) = C_{i+1}$, $f^{a_k}(C_i) = C_i$ and $f^{a_k}|_{C_i}$ 
is topological mixing. \vspace{2mm}

\noindent The proof of Theorem D is similar to the classical proof of the spectral decomposition theorem for an expansive homeomorphism. The main difference is that Theorem 5.1 below is employed to deal with periodic
points. \vspace{2mm}

\noindent Theorem 5.1: Let $f$ be as in Theorem D. Then there exists $\epsilon > 0$ such that if $p \in \text{Per}(f)$, then $W^s_{\epsilon}(p,f)\subset W^s(p,f)$ and $W^u_{\epsilon}(p,f) \subset W^u(p,f)$. 
 
\subsection{Avalanche dynamics} 
This article is joint work with Ana Rodrigues and appeared in {\em Stochastics and\ Dynamics}\ (2016), \cite{MR3265617}. Denker's interest in  biological models began around 2004 in a collaboration  with the physicist M.\ Herrmann supervising the PhD thesis of Anna Levina \cite{levina} on the EHE-model  for neuronal avalanches.

The avalanche map studied by Denker and Rodrigues in this article, modeling avalanches in neural dynamics, is a new type of dynamical system motivated by the product of circle rotations as considered in Levina's thesis \cite{levina} in the metric space context. For simplicity we use this model here, but go with the notations of the paper itself. To start, let $X:=\mathbb{T}^N=\mathbb{T}_1\times \cdots \times \mathbb{T}_N$, where each component is a flat circle, i.e., a torus $\mathbb{T}=\mathbb{R}/\mathbb{Z}$, consider component-wise rotations  $S:=S_{\delta_1} \times \cdots \times S_{\delta_N}$, 
where $S_{\delta_i}(x)=x+\delta_i \pmod1$,  that are reasonably small compared to $N$, i.e. $N\delta_i<1$. 
The number of factors $N$ resembles the number of neurons in some neuronal network.  Given open sets, think intervals, in each component $U_i \subset S_i$ one defines the \emph{avalanche transformation}  $R_A$ as follows. First  
let $0\leq A(x,k) \leq N$ be the number of components $x_i$ whose $k$ orbit with respect to $R_{\delta_i}$ meets the respective open set $U_i$. 
Note, that the sequence $\{A(x,k)\}_k$ is non-decreasing and define the \emph{avalanche size} to be $A(x)=\max \{A(x,k): \text{ for all } j \leq k \ A(x,j) \geq j \}$. Now define the avalanche transformation as $R_A(x)=R^{A(x)+1}(x)$. The more general metric analog in the paper defines the avalanche size function $A$ the same way. Because of its coordinate-wise definition, combinatorics applies and determines the properties of the avalanche transformation, starting with a description of sets where $A(x)$ is constant. Once this is completed, 
Denker and Rodrigues are able to determine the distribution of the avalanche size, i.e. determine $m(\{x \in X:\ A(x)=a\})$ depending on suitable partitions of $\{1,...,N\}$ (Theorem 3.6).  Here $m$ is the product measure. 

In section 4 the avalanche transformation is again assumed to be invertible. In order to describe the wandering set of 
the avalanche transformation Denker and Rodrigues first describe a (nontrivial) measurable set $B \subset X$ that is $S_A$ invariant, i.e., $S_A(B) =B$.  Then this set is used (Theorem 4.3), under the invertibility assumption of the component maps, to determine a Borel measurable set $F \subset X$ that is forward invariant with respect to 
$S_A$, so that $F^c$ is contained in the wandering set for $S_A$ and so that $S_A$ equals the induced transformation for $F$.  Moreover the relative probability on $F$, induced by the invariant product measure, defines an $S_A$ invariant measure on $F$. Based on Theorem 4.3, they can now give a criterion for ergodicity for the avalanche transformation $S_A$ that depends on the $L^1$ eigenvalues of the unitary transformation defined by the component maps $(X_i, S_i, m_i)$. Specifically, $S_A$ is 
ergodic if all component maps are ergodic and have no common eigenvalues other than 1 (Theorem 4.5). 
Now Denker and Rodrigues derive a general combinatorial formula for the distribution of the avalanche dynamics which is based on Cayley's theorem on the number of labeled trees with $N$ vertices. In fact, the paper starts with a novel proof of that theorem.  
As expected, Denker and Rodrigues conclude a central limit theorem for the avalanche transformation 
$S_A$ (Theorem 4.7).

\subsection{Perspective}

Denker’s work in dynamical systems is characterized by the systematic use of probabilistic and ergodic methods to study concrete classes of systems. A central theme is the analysis of dynamical systems through their statistical properties, most prominently central limit theorems, but also invariance principles, and other limit laws, often derived using techniques from probability theory.

At the same time, Denker contributed to the structural understanding of dynamical systems through concepts such as Markov fibred systems and their random counterparts. These provide a flexible framework that extends classical symbolic dynamics and allows the application of thermodynamic formalism in non-uniform and stochastic settings.

A further characteristic feature is the interplay between general theory and specific examples, including rational maps, skew product systems, and random dynamical systems. Through the use of transfer operators, Gibbs measures, and related tools, Denker’s work connects geometric, probabilistic, and dynamical aspects in a unified framework.

\section{Statistics}
\label{sec:stats}

Manfred Denker has made pioneering and lasting contributions to statistics, covering a wide range of topics, mostly in the area of non-parametric statistics. 
His interest in statistics was originally motivated by the idea to gain information on the properties of dynamical systems by observing a finite segment $X_k=f(T^k \omega)$, $1\leq k\leq n$, of observables of the orbit $(T^k \omega)_{k\geq 1}$. These could either be real life data or data generated by computer simulations from a given system, where it 
is impossible to make analytic calculations, e.g., about the dimension of an attractor or about properties of an invariant distribution such as symmetry. Making inference about a data generating process from a finite set of observations is a classical problem in statistics. 
In order to quantify the uncertainty of the statistical procedures, one needs to calculate the asymptotic distribution of estimators and test statistics. In connection with data from dynamical systems  the major challenge lies in the fact that the observations are not i.i.d.\ which is the standard assumption in classical large sample theory. Manfred Denker has made seminal contributions to asymptotic distribution theory in nonparametric statistics for dependent data. Already in his early papers, e.g. in the programmatic article ``Statistical Decision Procedures and Ergodic Theory" (1981) \cite{MR730767}, Manfred Denker emphasized that the existing results that focus on mixing processes are not suitable in the setting of dynamical systems, where the best one can hope for are functionals of absolutely regular processes, a concept known from ergodic theory. 

While most of his papers focus on weakly dependent data, Manfred Denker has also made profound contributions in the classical setting of i.i.d.\ data, e.g., in connection with the law of the iterated logarithm and strong approximation results for U-statistics,  and the central limit theorem for empirical $U$-processes. Many of his innovative approaches to large sample theory have found their way into his 1985 monograph ``Asymptotic Distribution Theory in Nonparametric Statistics" \cite{MR889896}.

\subsection{Rank Statistics} Most early work in statistics centered on Gaussian data, where linear procedures are optimal. It was not until the 1940s that statisticians became aware of the drawbacks of Gaussian based procedures, especially their sometimes dramatic loss of efficiency in the presence of non-Gaussian data. Suitable nonparametric procedures  have better efficiency for non-Gaussian data, while their loss of efficiency in the Gaussian case is limited. The first such class of nonparametric procedures are rank tests, where the basic idea is to replace the original data by their ranks in the sample. In the simplest one sample setting with continuous data $X_1,\ldots,X_n$, the smallest observation gets rank $1$, the second smallest rank $2$, and so on, until the largest observation gets rank $n$. Formally, the rank of the $i$-th observation can be defined by 
\[
 R_i= \# \{1\leq j\leq n: X_j\leq X_i  \}. 
\]
The functions of the ranks are called rank statistics, and the resulting tests are called rank tests. 

The best known examples of rank tests are Wilcoxon's two-sample rank test, Wilcoxon's signed rank test, and the Kruskal-Wallis test for homogeneity in the $k$-sample problem 
\cite{MR395032}.
Wilcoxon's two-sample rank test is designed for testing equality of the distributions $F$ and $G$ when $X_1,\ldots,X_{n_1}$ are $F$-distributed, and $X_{n_1+1},\ldots,X_{n_1+n_2}$ are $G$-distributed, against the alternative that $G$ is stochastically larger than $F$. The Wilcoxon test rejects the hypothesis that $F$ is equal to $G$ for large values of the rank sum
\[
 W=\sum_{i=n_1+1}^{n_1+n_2} R_i.
\]
Wilcoxon's two-sample rank test is a nonparametric competitor to the two-sample Student $t$-test. It is known to have an efficiency of $3/\pi\approx 95\%$ for Gaussian data, i.e., that it requires roughly $5\%$ more observations to have the same power as the $t$-test, while it sometimes dramatically outperforms the $t$-test for non-Gaussian data. 

Wilcoxon's two-sample rank test statistic is a special example of the class of so-called simple linear rank statistics, which are defined as
\[
  S= \sum_{i=1}^n c_n(i) a_n(R_i),
\]
where the $(c_n(i))_{1\leq i\leq n}$ are the regression constants and where $(a_n(i))_{1\leq i\leq n}$ are the scores. Under some assumptions on the scores and the regression constants, Wald and Wolfowitz (1944, \cite{MR11424}) established asymptotic normality of simple  linear rank statistics in case the data are independent and identically distributed. Chernoff and Savage (1958, \cite{MR100322}) determined the asymptotic distribution for independent, but not necessarily identically distributed data, thus making it possible to calculate the asymptotic relative efficiency of rank tests. 

Chernoff and Savage employ an empirical process representation of simple  linear rank statistics. Note that the ranks can be represented via the empirical distribution function $F_n:\R\rightarrow [0,1]$, which is defined by $F_n(x)=\frac{1}{n}\#\{1\leq i\leq n: X_i\leq x  \}$, so that one obtains
\[
  R_i=n\, F_n(X_i). 
\]
 It is usually assumed that the scores are given by a score function $h:[0,1]\rightarrow \R$ in such a way that 
$a_n(i)=h(\frac{i}{n})$. With these definitions, the simple  linear rank statistic $S$ can be expressed as follows
\[
S=\sum_{i=1}^n c_n(i) h(\frac{R_i}{n})  = \sum_{i=1}^n c_n(i) h(F_n(X_i)) 
=\int h(F_n(x)) dH_n(x),
\]
where $F_n$ is the empirical distribution function of the data $X_1,\ldots,X_n$, and where $H_n$ is the weighted empirical distribution function assigning weights $c_n(1),\ldots, c_n(n)$ to the data $X_1,\ldots,X_n$, i.e.,
$H_n(x)=\sum_{i=1}^n c_n(i) 1_{\{X_i\leq x\}}$.

In a series of papers beginning in the late 1970s,  Manfred Denker initiated a systematic study of the asymptotic distribution of simple linear rank statistics for weakly dependent data, building forth on the above mentioned empirical process representation of rank statistics. The key approach exploited by Denker was  to view the linear rank statistic as a linear function of the score function $h$ and to analyze the sequence of operators 
\[
  h\mapsto S_n(h): = \frac{1}{\sqrt{n}} \Big( \int h(F_n(x)) dH_n(x) - \int h(F(x)) dH(x) \Big) 
\]
as maps from a suitable normed linear space of score functions into the space of probability distributions, equipped with the Wasserstein distance $d_2$, defined by
\[
  d_2(P,Q)= \inf_{X\sim P,\, Y\sim Q} \sqrt{E(|X-Y|^2)}.
\]
Under suitable assumptions, this sequence of operators is uniformly bounded, and thus one can conclude that the set of score functions for which asymptotic normality holds is closed. It then remains to establish asymptotic normality for a dense set of simple score functions, which can usually be achieved by applying known central limit theorems for partial sums. Following this approach, Denker and R\"osler (1985) \cite{MR792458}
were able to prove asymptotic normality of simple linear rank statistics under conditions weaker than those in any previously known results. Brunner and Denker (1994) \cite{MR1309630} used the same approach to show asymptotic normality of rank tests in factorial designs. 

Both the collaboration with Uwe R\"osler and with Edgar Brunner took place locally  in G\"ottingen. During the  1970s and 1980s Uwe R\"osler was an undergraduate, later a Ph.D.\ student, and finally an assistant professor in the Institute of Mathematical Stochastics, before moving on to a permanent position as professor at the University of Kiel.  Edgar Brunner held the chair in Medical Statistics in the G\"ottingen Medical School from 1976 until his retirement in 2009. Because of his strong mathematical background, Brunner connected easily with his colleagues in the Math Department, in the first place with Manfred Denker. Their cooperation spanned more than 30 years, during which Brunner and Denker co-supervised many Master's and Ph.D.\ theses. Many of their research problems arose from Brunner's statistical consulting with colleagues in the medical school.

\subsection{Logarithmic quantile estimation} In a series of papers beginning in the early 2000s, Denker proposed the use of the pathwise central limit theorem as an alternative to the bootstrap for determining confidence intervals and critical values for tests. 

To explain the main idea in a simple setting, take the arithmetic mean $\bar{X}_n=\frac{1}{n} \sum_{i=1}^n X_i$ as an estimator for the expected value $\mu=\E(X)$. If one knew the distribution $G_n$ of $\sqrt{n}(\bar{X}_n-\mu)$, one could construct a $(1-\alpha)$-confidence interval for $\mu$ as
\[
  \big[\bar{X}_n-\frac{q_{n,1-\alpha/2}}{\sqrt{n}}, \bar{X}_n- \frac{q_{n,\alpha/2}}{\sqrt{n}} \big],
\] 
where $q_{n,p}$ denotes the $p$-th quantile of $G_n$, i.e., $G_n(q_{n,p})=p$. Except for the case when the random variables $X_i$ are i.i.d.\ Gaussian with known variance, the distribution of $G_n$ is generally unknown. 
A popular method for the calculation of approximate quantiles of $G_n$ is the bootstrap, where one estimates the distribution $G_n$ by resampling the statistic $\sqrt{n}(\bar{X}-\mu)$ from the empirical distribution of the data $X_1,\ldots,X_n$. 

As an alternative to the bootstrap, Manfred Denker proposed the use of the pathwise central limit theorem for estimation of the quantiles of $G_n$. The pathwise central limit theorem, also known as almost sure central limit theorem, was discovered in the late 1980s independently by Brosamler (1988) \cite{MR957261}, Fisher (1987) \cite{MR877784}, and Schatte (1988) \cite{MR968997}.
For i.i.d.\ random variables with finite variance $\sigma^2=\Var(X_1)$, the pathwise CLT states that almost surely, as $N\longrightarrow \infty$,
\[
  \widehat{G}_n(x):=\frac{1}{C_n} \sum_{k=1}^n \frac{1}{k} 1_{\{ \sqrt{k}(\bar{X}_k-\mu) \leq x  \}} \longrightarrow G(x),
\]
where $G$ denotes the distribution function of the $N(0,\sigma^2)$-distribution, and where $C_n=\sum_{k=1}^n \frac{1}{k} \sim \log n$. 
Initially, the pathwise CLT came as a big surprise, as it states that the CLT, which is a distributional law, can be observed along a single realization of the process $(X_i)_{i\geq 1}$. Berkes and Dehling (1993) \cite{MR1235433} showed that pathwise versions of distributional limit theorems also hold for convergence to a stable distribution.

Denker proposed to estimate the quantiles of $G$ by the corresponding quantiles of $\widehat{G}_n$. 
This idea was worked out in detail, and applied to the nonparametric Behrens-Fisher problem, in the 2005 
G\"ottingen Ph.D.\ thesis by Thangavelu, which was jointly supervised by Edgar Brunner and Manfred Denker. In a joint 2015 paper with his former Ph.D.\ student Lucia Tabacu \cite{MR3275537}, Denker established this method, called logarithmic quantile estimation, for rank statistics. The main ingredient in their paper is the proof of the almost sure central limit theorem for simple linear rank statistics.

\subsection{U-Statistics}
In a series of papers beginning in the early 1980s, and  motivated by his research on the statistical analysis of dynamical systems,  Denker investigated the large sample behavior of statistical functionals, i.e., of measurable  functions $T_n=t_n(X_1,\ldots,X_n)$ of the data $X_1,\ldots, X_n$. A famous result by Richard von Mises (1947) \cite{MR22330}
states that under very general conditions such functionals can be expanded into a series of von Mises statistics
\[
  V_n(h)=\frac{1}{n^m} \sum_{1\leq i_1 , \ldots , i_m\leq n} h(X_{i_1},\ldots,X_{i_m}),
\]
where the so-called kernel $h(x_1,\ldots,x_m)=h^{(m)}(x_1,\ldots,x_m)$ is a measurable function that is invariant under permutations of its $m$ arguments.  
The large sample behavior of von Mises statistics is closely related to that of the associated U-statistics, where the diagonal terms are deleted from the summation, 
\[
  U_n(h)= \frac{1}{\binom{n}{m}} \sum_{1\leq i_1< \ldots <i_m \leq n} h(X_{i_1},\ldots,X_{i_m}).
 \]
U-statistics were introduced independently by Halmos (1946)  \cite{MR15746} 
and Hoeffding (1948)  \cite{MR26294}
in connection with the theory of unbiased estimation. Note that in the case of i.i.d.\ data, $U_n(h)$ is an unbiased estimator of the parameter $\theta=Eh(X_1,\ldots,X_m)$.
 For ease of presentation, we will focus on the bivariate case, i.e.,  when $m=2$. The fundamental tool in the analysis of a U-statistic is the Hoeffding decomposition of the kernel, given by $h(x,y)=\theta + h_1(x)+h_1(y) +h_2(x,y)$, where the terms on the right-hand side are defined as
\begin{align*}
\theta&= E(h(X_1,\tilde{X})) \\
h_1(x)&=E(h(x,X_1))-\theta\\
h_2(x,y)&=h(x,y)-\theta -h_1(x)-h_1(y),
\end{align*}
where $\tilde{X}$ is an independent copy of $X_1$. 
The decomposition of the kernel gives rise to the Hoeffding decomposition of the U-statistic
\begin{align*}
U_n(h)=\theta + \frac{2}{n}\sum_{i=1}^n h_1(X_i) + U_n(h_2).
\end{align*}
The asymptotic behavior of $U_n(h)$ depends crucially on whether the kernel is degenerate, which is the case when $h_1(x)=0$ for all $x$, or whether it is non-degenerate. In the non-degenerate case, the asymptotic behavior of the U-statistic is governed by the first order term in the Hoeffding decomposition $\sum_{i=1}^n h_1(X_i)$, and thus one obtains the same results as for partial sum processes. In order to show that the first order term is dominating, one needs to control the second order term $U_n(h_2)$. For i.i.d.\ data, this can be done by standard variance estimates using the fact that $\sum_{1\leq i<j\leq n} h_2(X_i,X_j)$ is a martingale. For dependent data, the control of the second order term requires a more careful analysis, usually employing correlation inequalities for the summands $h_2(X_i,X_j)$. 

In his research, Denker made lasting contributions  to the large sample behavior of U-statistics, both for i.i.d.\ as well as for dependent data, and both for degenerate and non-degenerate kernels. In the following, we discuss some of his main achievements in more detail.

\subsubsection{Almost sure approximation of degenerate U-statistics} In two papers \cite{MR758070, MR835162} written in the early 1980s, Manfred Denker together with co-authors Herold Dehling and Walter Philipp
obtained pioneering results on the fluctuation behavior of U-statistics of i.i.d.\ data with a degenerate kernel. In the first of these papers, devoted to the bivariate case, they proved an almost sure approximation of a degenerate U-statistic with square integrable kernel by a stochastic integral with respect to a Kiefer-M\"uller process.  More specifically, they could show that, after possibly enlargening the probability space, there exists a Kiefer-M\"uller process $K(x,t)$ such that 
\[
  \sum_{1\leq i <j\leq n} h(X_i,X_j) - \iint h(x,y) K(dx,n) K(dy,n) = o(n\, \log \log n).
\]
This theorem constituted the first almost sure approximation result beyond the domain of partial sums. In addition, the authors could establish the  law of the iterated logarithm for degenerate U-statistics, stating that there exists a constant $C(h)\in (0,\infty)$ such that 
\[
  \limsup_{n\rightarrow \infty} \frac{1}{n \log\log n} \sum_{1\leq i <j\leq n} h(X_i,X_j) = C(h).
\]
In fact, they could identify $C(h)$ as the largest eigenvalue of the integral operator with kernel $h$. 
The proof of the Dehling-Denker-Philipp (1984)
almost sure approximation theorem combined then modern ideas from Probability in Banach Spaces, especially the Law of the Iterated Logarithm in Hilbert spaces, together with martingale limit theorems. In a subsequent paper, published in 1986, the same authors were  able to extend their results to degenerate U-statistics of arbitrary order. This extension required subtle new results on the law of the iterated logarithm for Hilbert-space valued martingales. 
In a later paper, Dehling (1989)  \cite{MR991944}
established  a functional version of the law of the iterated logarithm, identifying the almost sure limit set of the sequence of functions $ t\mapsto \frac{1}{n \log\log n} \sum_{1\leq i <j\leq [nt]} h(X_i,X_j)$, $t\in [0,1]$. 
In addition,  Dehling (1989) \cite{MR991944}
was able to identify the $\limsup$ in the law of the iterated logarithm for U-statistics of arbitrary degree $m$.  

\subsubsection{U-statistics of dependent data}
In two papers, written in the early 1980s, and both coauthored with Gerhard Keller, Manfred Denker investigated the asymptotic distribution of U-statistics of dependent data in the non-degenerate case. At the time of writing these papers, Gerhard Keller was a postdoc at G\"ottingen, where he had arrived in 1979 after having obtained his Ph.D.\ at Rennes under the supervision of Mike Keane. He finally became professor in Erlangen. 
In their first paper, Denker and Keller (1983) \cite{MR717756} obtained mostly theoretical results such as central limit theorems, Berry-Esseen bounds, functional central limit theorems and almost sure approximations by a Wiener process under weak assumptions on the moments of the kernel $h$ and on the mixing coefficients. A typical condition imposed in this paper is that the process $(X_j)_{j\geq 1}$ should be $\beta$-mixing with mixing coefficients satisfying $\sum_{n=1}^\infty \beta(n)^{\delta/(2+\delta)}<\infty$ and that 
\[
  \sup_{1\leq i_1<\ldots <i_m } \E(|h(X_{i_1},\ldots, X_{i_m})|^{2+\delta})<\infty.
\]
In their second paper, Denker and Keller (1986) \cite{MR854400} studied applications of U-statistics to the statistical analysis of dynamical systems. 
One of the motivating examples was the Grassberger-Procaccia estimator for the correlation dimension of the attractor of a dynamical system which is defined as the exponent $d$ in the asymptotic expansion $C(r):=P(|X-X^\prime |\leq r)\sim r^d$, where $X$ and $X^\prime$ are independent observations from the invariant distribution. The sample analogue of the correlation integral $C(r)$ is the 
 U-statistic
\[
  U_n(r)=\frac{1}{\binom{n}{2}} \sum_{1\leq i<j\leq n} 1_{\{|X_i-X_j|\leq r\}}.
\]
Using the relation $\log C(r) \approx d \log r$, it is natural to estimate the correlation dimension $r$ by linear least squares methods applied to the points $(\log r_1,\log U_n(r_1)), \ldots$,  $(\log r_k, \log U_n(r_k))$.  
Denker and Keller (1986)  establish joint asymptotic normality of $U_n(r_1),\ldots, U_n(r_k) $ in case the underlying data are functionals of an absolutely regular process, i.e., that $X_j$ can be represented as 
\[
  X_j=f(Z_j, Z_{j+1}, Z_{j+2},\ldots),
\]
for some Lipschitz-continous function $f:\mathcal{X}^\N \rightarrow \R$ and some absolutely regular process $(Z_j)_{j \geq 1}$. This condition is an important weakening of any of the classical mixing assumptions, and is satisfied by a large class of dynamical systems such as expanding maps of the unit interval; see, e.g. Hofbauer and Keller (1982) \cite{MR656227}
and Burton, Borovkova and Dehling (2001) \cite{MR1851171}.
As an application, they could show asymptotic normality of the least squares estimator of $d$. More generally, Denker and Keller (1986) could show asymptotic normality of U-statistics with general kernels $h(x,y)$, under some continuity assumptions on the kernel. 

\subsubsection{Empirical U-processes}
The sample correlation integral $U_n(r)$ is on one hand a U-statistic, and on the other hand also an empirical process of the sample distances $|X_i-X_j|$, $1\leq i<j\leq n$. Replacing the sample distances by arbitrary functions $g(x,y)$, one obtains the so-called empirical distribution function of U-statistics type, defined as
\[
  U_n(x)=\frac{1}{\binom{n}{2}} \# \{1\leq i< j\leq n: g(X_i,X_j)\leq x\} =\frac{1}{\binom{n}{2}} \sum_{1\leq i<j\leq n} 
  1_{\{g(X_i,X_j) \leq x\}}.
\]
In a joint paper with Herold Dehling and Walter Philipp, Manfred Denker (1987)\cite{MR891707}
could establish convergence in distribution of the so-called empirical U-process 
\[
 \sqrt{n} (U_n(x)-P(g(X_1,X_2)\leq x)) 
\]
to a mean-zero Gaussian process, thus generalizing Donsker's classical empirical process invariance principle to empirical distributions of U-statistics type. In their paper, Dehling, Denker and Philipp (1987) assumed that the underlying data $(X_i)_{i\geq 1}$ are independent and identically distributed. Later, Borovkova, Burton and Dehling (2001) could extend these results to dependent data, especially to functionals of absolutely regular processes, thus covering a wide range of examples from dynamical systems. 

\subsubsection{V-statistics of measure preserving transformations}
In a seminal paper, co-authored with his long-term collaborator Mikhail Gordin, Manfred Denker  investigated the asymptotic behavior of von Mises statistics of a measure preserving transformation $T$ on a probability space $(\mathcal{X},\mathcal{F},P)$, defined by
\[
  x\mapsto \frac{1}{C_n} \sum_{0\leq i_1,\ldots, i_m \leq n} f(T^{i_1}x,\ldots,T^{i_m}x),
\]
where $f:\mathcal{X}^d \rightarrow \R$ is a measurable function, and where $C_n$ are norming constants.  In the first place, Denker and Gordin (2014) \cite{MR3256808} raise the fundamental question about an appropriate function space for the kernels $f$, such that the above V-statistic is well-defined. The naive choice, to view $f$ as a representative of an equivalence class in the space $L_p(\mathcal{X}^m,\mathcal{F}^{\otimes m}, P^m)$ does not work, because this does not provide a proper definition of  $f$ on any of the diagonals that are null sets with respect to the product measure. Denker and Gordin (2014) show that the projective tensor product of $L_p(\mathcal{X},\mathcal{F},P)$ solves this problem. For such kernels, they establish the individual ergodic theorem for V-statistics, providing a partial answer to a question that had been raised by Aaronson, Burton, Dehling, Gilat, Hill, and Weiss (1996) \cite{MR1363941}, who had presented counterexamples to the ergodic theorem for U- and V-statistics, together with some sufficient conditions for the ergodic theorem to hold. In addition, Denker and Gordin (2014) establish distributional convergence of V-statistics of measure-preserving transformations using martingale techniques. 

\subsection{Parametric Inference for Dynamical Systems} 
While most of Manfred Denker's work in statistics is devoted to nonparametric methods, there is a recent contribution to a parametric estimation and testing problem. In the paper ``Parametrized families of Gibbs measures and their statistical inference'', written jointly with Marc Ke\ss eb\"oh\-mer, Artur O.\ Lopes and  Silvia R.C.\ Lopes, and published in \emph{Stochastics and Dynamics} (2025, \cite{MR4901113}),  the authors  investigate  maximum likelihood estimation and hypothesis testing  for unknown parameters in a parametrized family of Gibbs measures.
For known H\"older continuous functions $f_0,\ldots,f_d$ on a subshift of finite type, the authors consider the family of potentials $(F_\theta)_{\theta \in \Theta}$, defined by
\[
  F_\theta=f_0+\sum_{i=1}^d \theta_i f_i,
\]
together with the corresponding Gibbs measures. 
In their paper, the authors show that the maximum likelihood estimator (MLE) $\widehat{\theta}_n$ for $\theta$ is $\sqrt{n}$-consistent, and they identify its asymptotic  distribution. More precisely, under the assumption that the functions $f_1,\ldots, f_d$ are linearly independent as cohomology classes, the authors show that 
$\sqrt{n}(\widehat{\theta}_n -\theta)$ converges in distribution to $G^{-1}(N)\, N^t$. Here $G(N)=N^t\, N -\Sigma^{\mu_\theta}$, and $N$ is a centered multivariate normally distributed random vector with covariance matrix $\Sigma^{\mu_\theta}$, which is the limit covariance matrix of the functions $f_1,\ldots, f_d$ under the shift invariant Gibbs measure $\mu_\theta$ corresponding to the potential $F_\theta$. Moreover, the authors show that the MLE is asymptotically efficient in the sense that the mean squared error attains the Cram\'er-Rao lower bound. As an application of their results about the MLE, the authors  study the generalized likelhood ratio test for testing hypotheses about the unknown parameter $\theta$.

\subsection{Perspective}
Denker’s work in statistics focuses on the extension of classical asymptotic methods to dependent data, particularly in the context of dynamical systems. A central theme is the analysis of statistical procedures such as rank statistics, $U$--statistics, and empirical processes under weak dependence assumptions. By combining tools from probability theory, ergodic theory, and empirical process methods, Denker established limit theorems for statistical functionals that apply far beyond the independent setting. These contributions provide a framework for statistical inference in dynamical and other dependent environments, and illustrate the interplay between statistics, probability, and dynamical systems.

\section{Computer-Aided Results}
\label{sec:computer-aided}

Although Manfred Denker is widely recognized for his foundational contributions to ergodic theory, probability, and dynamical systems, an important aspect of his work is the development of mathematically rigorous methods that connect theoretical analysis with computation and data. These contributions arise naturally from his broader research program, particularly in the context of statistical inference for dynamical systems and the study of limit theorems beyond classical settings.

Rather than serving as purely heuristic tools, computational and algorithmic methods in Denker’s work are closely tied to rigorous probabilistic and ergodic arguments. They provide a bridge between abstract theory and practical questions involving finite data, numerical simulation, and statistical estimation.

\subsection{Statistical inference for chaotic dynamics}

Denker's early engagement with computer--aided methods appears already in the
mid--1980s, most notably in the paper
``Rigorous statistical procedures for data from dynamical systems" (\textit{Journal of Statistical Physics}, 1986)
\cite{MR854400}, coauthored with Gerhard Keller.  This work addressed a problem of central relevance to both
theory and applications: how to infer invariant statistical properties of a
dynamical system from finite observational data when the underlying dynamics are
deterministic but chaotic. In particular, they show that estimators such as the Grassberger--Procaccia correlation integral can be treated as $U$--statistics, allowing one to derive asymptotic distributions and quantify statistical errors rigorously. 

The conceptual framework of this paper is closely connected to Denker's broader
work on mixing properties, weak Bernoulli systems, and invariance principles for
dependent sequences developed during the late 1970s and early 1980s.  By
exploiting probabilistic representations of chaotic trajectories, the authors
demonstrated that classical statistical tools can be used in a mathematically
sound manner for deterministic systems with sufficient stochastic structure. This approach is closely connected to Denker's broader work on weak dependence and invariance principles, where deterministic dynamical systems are modeled as stochastic processes with controlled dependence, allowing statistical procedures to be justified rigorously.

A key application treated in the paper is the rigorous estimation of geometric
characteristics of attractors, including the Grassberger--Procaccia correlation
dimension.  While this dimension had already gained popularity in the applied
literature through numerical experimentation, Denker and Keller (1986)  embedded
its estimation into the theory of $U$--statistics, thereby providing a precise
statistical interpretation. For example, the correlation integral can be written
in the form
\[
U_n(r)=\frac{1}{\binom{n}{2}} \sum_{1\leq i<j\leq n}
\mathbf{1}_{\{|X_i-X_j|\le r\}},
\]
which places the problem into the framework of $U$--statistics and allows one to
derive asymptotic normality and variance estimates under suitable weak
dependence assumptions. In this direction,  the major breakthrough of this paper is the proof of asymptotic normality of a large class of $U$-statistics when the underlying data might not be mixing, but can be represented as functionals of mixing processes, a structure  frequently found in chaotic dynamical systems. 

 This work set an early
standard for how computer--generated data from dynamical systems should be
analyzed when mathematical rigor is required. Together with the companion paper ``On U-statistics and von Mises' statistics for weakly dependent processes'' (\emph{Zeitschift f\"ur Wahrscheinlichkeitstheorie und verwandte Gebiete} 1983) \cite{MR717756}, also written jointly with Gerhard Keller,  this work has had an enormous impact on the development of nonparametric statistics of dependent data---as is also witnessed from the hundreds of citations that the papers have received, and continue to receive even 40 years after their  publication. 

\subsection{Stable resampling and non--Gaussian limits}

A second, and longer--running strand of Denker's work, spanning from the mid 1980s
to the early 2020s, concerns resampling methods based on stable laws and their
application to statistics with heavy--tailed behavior.  These ideas first
appeared in the paper ``{Resampling $U$--statistics using $p$--stable laws}''
(1990) \cite{MR1062544}, coauthored with Herold Dehling and Wojbor Woyczy\'nski (1943--2021), and published in the \textit{Journal of Multivariate Analysis}. This paper was written in the summer of 1985 when Dehling was postdoc at the Institute of Mathematical Statistics in G\"ottingen, and Woyczy\'nski visited the Institute. Denker and Woyczy\'nski had known each other since the 1975 Oberwolfach {\em Probability in Banach Spaces} meeting, which marked the beginning of their life-long friendship. In 1998 they published  the textbook {\em Introductory Statistics and Random Phenomena -- Uncertainly, Complexity and Chaotic Behavior in Engineering and Science} \cite{MR1643201}, which grew out of an NSF-funded initiative to restructure the engineering curriculum by incorporating the latest technological innovations.

In \cite{MR1062544}, resampled $U$-statistics are constructed whose weak limits are
expressed as multiple stochastic integrals with respect to $p$--stable motions.
This approach replaces classical Gaussian limits by stable laws, thereby extending
invariance principles to settings without finite second moments. More precisely,
the resampled U-statistics are defined as multilinear forms 
\[
\sum_{1\leq i_1<\ldots<i_m\leq n}  h(X_{i_1},\ldots,X_{i_m}) \, Y_{i_1} \cdot \ldots \cdot Y_{i_m},
\]
where $(Y_i)$ belongs to the domain of attraction of a $p$--stable law, leading to
limits described by multiple stable stochastic integrals. 
At the same time, this provides
a concrete resampling framework that overcomes limitations of standard bootstrap
methods in the presence of heavy tails, and yields consistent estimators even
when variance-based normalization fails.

This contribution is deeply connected to Denker's earlier work on $U$--statistics,
von Mises functionals, and non--Gaussian limit theorems from the early 1980s. 
The proof relies on empirical process techniques, establishing convergence of the resampled empirical process
$F_n(x)=\sum_{i=1}^n Y_i\, 1_{\{X_i \leq x \}}$ to a $p$-stable motion. 
 In
contrast to purely theoretical results, however, the resampling procedure is
designed to be directly implementable, with convergence results strong enough to
support numerical approximations and simulations.

More than thirty years later, these ideas were developed further in
``{Estimation by stable motions and its applications}" (2023) \cite{MR4660830}.
This paper develops explicit algorithms based on almost sure limit theorems for
stochastic integrals with respect to stable motions. More precisely, the method
constructs almost sure approximations to the limiting distribution using
logarithmic averages of the form
\[
\frac{1}{\log n}\sum_{k=1}^n \frac{1}{k}\,
\mathbf{1}_{\{S_k \le x\}},
\]
which converge almost surely to the limiting distribution. This allows the
estimation of quantiles and confidence intervals even in situations with heavy-tailed
data where classical Gaussian approximations are not applicable.

Taken together, these works illustrate a rare continuity between abstract
probability theory and computational methodology: ideas introduced in the mid-1980s are revisited, refined, and transformed into practical statistical tools
in the 2020s.

\subsection{Algorithmic construction of conformal measures}

A further domain in which Denker combined rigorous theory with computation is the
study of conformal measures in complex dynamics.  Following his earlier
fundamental contributions to the existence and uniqueness of conformal measures
for rational maps in the 1980's and 1990's, Denker turned to the question of how such
measures can be accessed numerically.

This perspective culminated in the paper
``{Pseudorandom numbers for conformal measures}" (2009) \cite{MR2572997}.  The authors propose an algorithm for constructing pseudo--generic points for conformal and invariant measures associated with expanding maps and rational functions. The approach is based on discretized conformal equations combined with least--squares optimization and avoids direct spectral computations for
transfer operators. More precisely, the construction approximates solutions to the conformal relation
\[
m(T(A))=\int_A e^{\phi}\, dm,
\]
for suitable potentials $\phi$, which characterizes conformal measures in the dynamical system. The authors construct sequences of points whose empirical distributions approximate the conformal measure, and establish convergence under suitable assumptions on the dynamics and the potential. From a dynamical point of view, this provides a numerical realization of conformal measures that bypasses direct spectral analysis of the associated transfer operator.

The algorithm is carefully analyzed with respect to convergence, numerical
accuracy, and robustness, and it is implemented explicitly for Julia sets of
hyperbolic quadratic polynomials.  This work demonstrates how theoretical
concepts such as conformal measures, invariant densities, and Hausdorff
dimension which are central to Denker’s earlier research, can be realized through controlled
numerical procedures.

\subsection{Visualization and computational exploration}

Denker’s interest in computation and visualization is also reflected in the
contribution ``{Polynomial skew products}" (2001) \cite{MR1850306}. This paper studies dynamical systems of the form
\[
F(z,w) = (p(z), q(z,w)),
\]
where $p$ is a polynomial in one complex variable and $q$ is a polynomial in two variables, leading to a skew product structure over the base dynamics defined
by $p$. Such systems naturally give rise to invariant sets whose geometry is determined by the interaction between the base dynamics and the fiber maps.

In particular, the associated Julia sets can be described in terms of fiberwise Julia sets over points in the base, resulting in fractal structures of high
complexity. These sets exhibit behavior that is difficult to access analytically, but can be explored through numerical and visual methods, which reveal their
structure and dynamical features.

Although partly expository, the paper illustrates how visualization and computational experimentation can play an essential role in understanding higher-dimensional dynamical systems, where explicit analytic descriptions are often not available.

\subsection{Perspective}
Denker’s work in this area reflects a consistent effort to extend the scope of rigorous mathematical analysis to settings involving finite data and numerical approximation. Rather than developing computation as an independent direction, he integrates algorithmic and statistical methods into the broader framework of ergodic theory and probability.

In this way, his contributions illustrate how abstract concepts, such as invariance principles, conformal measures, and stochastic limit laws, can be adapted to yield quantitative and computationally accessible results. This interaction between theory and computation remains an important direction in the study of dynamical systems and stochastic processes.

\bibliographystyle{amsplain}

\providecommand{\bysame}{\leavevmode\hbox to3em{\hrulefill}\thinspace}
\providecommand{\MR}{\relax\ifhmode\unskip\space\fi MR }
\providecommand{\MRhref}[2]{%
  \href{http://www.ams.org/mathscinet-getitem?mr=#1}{#2}
}
\providecommand{\href}[2]{#2}

\end{document}